# Prescribing a fourth order conformal invariant on the standard sphere, Part II: blow up analysis and applications

Zindine Djadli, Andrea Malchiodi and Mohameden Ould Ahmedou

**Abstract**

In this paper we perform a fine blow up analysis for a fourth order elliptic equation involving critical Sobolev exponent, related to the prescription of some conformal invariant on the standard sphere $(\mathbb{S}^n, h)$. We derive from this analysis some a priori estimates in dimension 5 and 6. On $\mathbb{S}^5$ these a priori estimates, combined with the perturbation result in the first part of the present work, allow us to obtain some existence result using a continuity method. On $\mathbb{S}^6$ we prove the existence of at least one solution when an index formula associated to this conformal invariant is different from zero.

## 1 Introduction

Let $(M, g)$ be a smooth 4-dimensional Riemannian manifold, and consider the following fourth order operator, introduced by Paneitz [34]

$$P_g^4 \varphi = \Delta_g^2 \varphi - div_g \left( \frac{2}{3} Scal_g g - 2 Ric_g \right) d\varphi,$$

where $Scal_g$ and $Ric_g$ denote the scalar curvature and the Ricci curvature of $(M^4, g)$ respectively. Originally this operator was introduced for physical motivations and has many applications in mathematical physics, see [23], [35], [1] and [26]. Moreover the Paneitz operator enjoys many geometric properties, in particular it is conformally invariant, and can be seen as a natural extension of the Laplace-Beltrami operator on 2-manifolds. After the pioneering works by T. Branson [12], [13] and S.A. Chang and P. Yang [21], this operator has been the subject of many papers [14], [17], [18], [29], see also the surveys [15] and [22]. $P_g^4$ has been generalized to manifolds of dimension greater than 4 by T. Branson [13], and it is given for $n \geq 5$ by

$$P_g^n u = \Delta_g^2 u - div_g \left( a_n Scal_g g + b_n Ric_g \right) du + \frac{n-4}{2} Q_g^n u,$$

where

$$a_n = \frac{(n-2)^2 + 4}{2(n-1)(n-2)}, \quad b_n = -\frac{4}{n-2}$$

$$Q_g^n = \frac{1}{2(n-1)} \Delta_g Scal_g + \frac{n^3 - 4n^2 + 16n - 16}{8(n-1)^2(n-2)^2} Scal_g^2 - \frac{2}{(n-2)^2} |Ric_g|^2.$$

As for $P_g^4$, the operator $P_g^n$, $n \geq 5$, is conformally invariant; if $\tilde{g} = \varphi^{\frac{4}{n-4}} g$ is a conformal metric to $g$, then for all $\psi \in C^\infty(M)$ we have
$$P_g^n(\psi\varphi) = \varphi^{\frac{n+4}{n-4}} P_{\tilde{g}}^n(\psi).$$

In particular, taking $\psi = 1$, it is

(1) $$P_g^n(\varphi) = \frac{n-4}{2} Q_{\tilde{g}}^n \varphi^{\frac{n+4}{n-4}}.$$

In this paper we continue to study the problem of prescribing $Q$ on the standard sphere $(S^n, h)$, $n \geq 5$. By equation (1), the problem can be formulated as follows. Given $f \in C^2(\mathbb{S}^n)$, we look for solutions of

(P) $$P_h^n u = \frac{n-4}{2} f u^{\frac{n+4}{n-4}}, \quad u > 0, \quad \text{on } \mathbb{S}^n.$$

On the unit sphere $(\mathbb{S}^n, h)$, $n \geq 5$, the operator $P_h^n$ is coercive on the Sobolev space $H_2^2(\mathbb{S}^n)$, and has the expression
$$P_h^n u = \Delta_h^2 u + c_n \Delta_h u + d_n u,$$
where
$$c_n = \frac{1}{2}(n^2 - 2n - 4), \quad d_n = \frac{n-4}{16} n (n^2 - 4).$$

Problem (P) is the analogous, for Paneitz operator, of the so-called Scalar Curvature Problem, to which many works are devoted, see [4], [6], [2], [8], [9], [10], [11], [20], [16], [27], [32], [30], [33], [38], and the references therein. We also refer to the monograph [5].

Our aim is to give sufficient conditions on $f$ such that problem (P) possesses a solution. It is easy to see that a necessary condition for solving (P) is that $f$ has to be positive somewhere. Moreover, there are also obstructions of Kazdan-Warner type, see [24] and [40].

In the first part of the present work, [25], we established a perturbative result, see for a particular case Theorem 1.1 below. Before stating Theorem 1.1, we introduce some preliminary notation following [20]. For $P \in \mathbb{S}^n$ and $t \in [1, +\infty)$ we denote by $\varphi_{P,t}$ the conformal map on $\mathbb{S}^n$ defined as follows: using stereographic coordinates with projection through the point $P$, we set
$$\varphi_{P,t}(y) = t\, y.$$

Let also $G: B^{n+1} \to \mathbb{R}^{n+1}$ be given by
$$G\left(\frac{t-1}{t} P\right) = \frac{1}{\omega_n} \int_{\mathbb{S}^n} (f \circ \varphi_{P,t}(x)) \, x \, dv(h),$$
where $\omega_n$ denotes the volume on $\mathbb{S}^n$ with respect to its standard volume $dv(h)$.

**Theorem 1.1** *There exists $\varepsilon(n) > 0$, depending only on $n$, such that (P) admits a solution for all $f \in C^\infty(\mathbb{S}^n)$ satisfying*

(ND) $$\Delta_h f(q) \neq 0 \quad \text{whenever} \quad \nabla f(q) = 0,$$

*and*
$$\left\| f - \frac{n(n^2-4)}{8} \right\|_{L^\infty(\mathbb{S}^n)} \leq \varepsilon(n),$$
$$\deg\left(G, \{(P,t) \in \mathbb{S}^n \times (1; +\infty), t < t_0\}, 0\right) \neq 0 \quad \text{for } t_0 \text{ large enough.}$$

Our main goal in this second part is to perform a fine blow up analysis of equation $(P)$ and of its subcritical approximation. Then we take advantage of this study to derive some compactness and non-perturbative existence results for problem $(P)$ in lower dimensions. In order to describe the blow up analysis we need some preliminaries.

Let $\Omega \subseteq \mathbb{R}^n$ be a bounded smooth domain, $\tau_i \geq 0$ satisfy $\lim_i \tau_i = 0$, $q_i = \frac{n+4}{n-4} - \tau_i$ and $\{f_i\}_i \in C^1(\Omega)$ satisfy

(2) $$\frac{1}{A_1} \leq f_i(x) \leq A_1, \quad \text{for all } x \in \Omega, \text{ for all } i,$$

for some positive constant $A_1$. Consider the family of equations

(3) $$\Delta^2 u_i = \frac{n-4}{2} f_i(x) u_i^{q_i}, \quad u_i > 0, \quad \Delta u_i > 0, \quad x \in \Omega.$$

We recall that according to our notation

$$\Delta \varphi = -\sum_{i=1}^{n} \frac{\partial^2 \varphi}{\partial x_i^2}, \qquad \varphi \in C^2(\Omega).$$

Our aim is to describe the behavior of $u_i$ when $i$ tends to infinity. It is possible to prove, see Remark 2.8, that if $\{u_i\}_i$ remains bounded in $L^\infty_{loc}(\Omega)$, then for any $\alpha \in (0,1)$ $u_i \to u$ in $C^{4,\alpha}_{loc}(\Omega)$ along some subsequence. Otherwise, we say that $\{u_i\}_i$ *blows up*. In the following we adapt to this new framework the definition of isolated blow up points and isolated simple blow up points which were introduced by R. Schoen, [37] (see also Y.Y. Li, [30]). Let $B_r(x) = \{y \in \mathbb{R}^n : |y - x| < r\}$.

**Definition 1.2** *Suppose that $\{f_i\}_i$ satisfy (2), and $\{u_i\}_i$ satisfy (3). A point $\overline{y} \in \Omega$ is called a blow up point for $\{u_i\}_i$ if there exists $y_i \to \overline{y}$, such that $u_i(y_i) \to +\infty$.*

In the sequel, if $\overline{y}$ is a blow up point for $\{u_i\}_i$, writing $y_i \to \overline{y}$ we mean that, for all $i$, $y_i$ is a local maximum of $u_i$ and $u_i(y_i) \to +\infty$ as $i \to +\infty$.

**Definition 1.3** *Assume that $y_i \to \overline{y}$ is a blow up point for $\{u_i\}_i$. The point $\overline{y} \in \Omega$ is called an* isolated *blow up point if there exist $\overline{r} \in (0, dist(\overline{y}, \partial\Omega))$ and $\overline{C} > 0$ such that*

(4) $$u_i(y) \leq \overline{C} |y - y_i|^{-\frac{4}{q_i - 1}}, \quad \text{for all } y \in B_{\overline{r}}(y_i) \cap \Omega.$$

Let $y_i \to \overline{y}$ be an isolated blow up point for $\{u_i\}_i$, we define $\overline{u}_i(r)$ to be (here $|\partial B_r|$ is the $n-1$-dimensional volume of $\partial B_r$)

(5) $$\overline{u}_i(r) = \frac{1}{|\partial B_r|} \int_{\partial B_r(y_i)} u_i, \quad r \in (0, dist(y_i, \partial\Omega)),$$

and

$$\hat{u}_i(r) = r^{\frac{4}{q_i - 1}} \overline{u}_i(r), \quad r \in (0, dist(y_i, \partial\Omega)).$$

**Definition 1.4** *An isolated blow up point $\overline{y} \in \Omega$ for $\{u_i\}_i$ is called an* isolated simple *blow up point if there exist some $\varrho \in (0, \overline{r})$, independent of $i$, such that $\hat{u}_i(r)$ has precisely one critical point in $(0, \varrho)$ for large $i$.*

The study of equation (3) has the following motivation. Taking a subcritical approximation of $(P)$, we consider

(6) $$P_h^n v_i - \frac{n-4}{2} f_i(x) v_i^{q_i} = 0, \quad v_i > 0, \quad \text{on } \mathbb{S}^n;$$

$$q_i = \frac{n+4}{n-4} - \tau_i, \quad \tau_i \geq 0, \quad \tau_i \to 0.$$

After performing a stereographic projection $\pi : \mathbb{S}^n \to \mathbb{R}^n$ through the north pole on $\mathbb{S}^n$, equation (6) is transformed into

$$\Delta^2 u_i(y) = \frac{n-4}{2} f_i(y) H_i(y)^{\tau_i} u_i(y)^{q_i}, \quad u_i > 0, \quad \text{on } \mathbb{R}^n,$$

where

(7) $$u_i(y) = \left(\frac{2}{1+|y|^2}\right)^{\frac{n-4}{2}} v_i(x), \quad H_i(y) = \left(\frac{2}{1+|y|^2}\right)^{\frac{n-4}{2}}, \quad y = \pi(x).$$

One can check that, under the assumption $f > 0$, if $v_i$ satisfies (6) and $u_i$ is given by (7), then it must be $u_i > 0, \Delta u_i > 0$ on $\mathbb{R}^n$, so we are reduced to study equation (3). Dealing with the blow up phenomenon we will use the same terminology for both $u_i$ and $v_i$, taking into account the relation (7).
Our main result regarding the blow up analysis for equation (6) is the following.

**Theorem 1.5** *Let $n = 5, 6$, and assume that $\{f_i\}_i$ is uniformly bounded in $C^1(\mathbb{S}^n)$ and satisfy (2). If $n = 6$ we also assume that $\{f_i\}_i$ is uniformly bounded in $C^2(\mathbb{S}^n)$. Let $\{v_i\}_i$ be solutions of (6): then*

$$\|v_i\|_{H_2^2(\mathbb{S}^n)} \leq C,$$

*where $C$ is a fixed constant depending on $n, A_1$ and the $C^1$ bounds of $\{f_i\}_i$ (and also the $C^2$ bounds of $\{f_i\}_i$ if $n = 6$). Furthermore, after passing to a subsequence, either $\{v_i\}_i$ remains bounded in $L^\infty(\mathbb{S}^n)$, or $\{v_i\}_i$ has only isolated simple blow up points, and the distance between any two blow up points is bounded below by some fixed positive constant.*

A fundamental ingredient in the blow up analysis is a Harnack type inequality for the above fourth order operator, proved in Lemma 2.5.
The blow up analysis can be specialized further in the case $n = 5$, yielding to a concentration phenomenon at most at one point for sequences of solutions of (6) and also to a compactness result for solutions of equation $(P)$ under condition $(ND)$.

**Theorem 1.6** *Let $n = 5$, and assume that $\{f_i\}_i$ is uniformly bounded in $C^1$ norm and satisfy (2). Let $\{v_i\}_i$ be solutions of (6). Then, after passing to a subsequence, either $\{v_i\}_i$ is bounded in $L^\infty(\mathbb{S}^n)$ or $\{v_i\}_i$ has precisely one isolated simple blow up point.*

**Theorem 1.7** *Let $n = 5$ and $f \in C^2(\mathbb{S}^5)$ be a positive Morse function satisfying $(ND)$. Then for any $\varepsilon > 0$ and any $\alpha \in (0,1)$ there exist positive constants $C = C(f,\varepsilon)$ and $\tilde{C} = \tilde{C}(f,\varepsilon,\alpha)$ such that for all $\mu$ with $\varepsilon \leq \mu \leq 1$, any positive solution $v$ of $(P)$ with $f$ replaced by $f_\mu = \mu f + (1-\mu)\frac{105}{8}$ satisfies*

$$C^{-1} \leq v(q) \leq C, \quad \text{for all } q \in \mathbb{S}^n; \quad \|v\|_{C^{4,\alpha}(\mathbb{S}^5)} \leq \tilde{C}.$$

Using Theorem 1.7, we derive the following existence result.

**Theorem 1.8** *Under the assumptions of Theorem 1.7, assume that the following condition holds,*

(8) $$\sum_{x \in Crit(f), \Delta_h f > 0} (-1)^{m(f,x)} \neq -1,$$

*where $Crit(f) = \{x \in \mathbb{S}^n | \nabla f(x) = 0\}$ and $m(f,x)$ denotes the Morse index of $f$ at $x$. Then equation $(P)$ has at least one solution, and the set of solutions of $(P)$ is compact in $C^{4,\alpha}(\mathbb{S}^5)$.*

The above Theorem is proved using a topological degree argument, following [16]. Theorem 1.1 provides the initial step of a continuity argument, while the a priori estimates in Theorem 1.7 are used to verify the invariance of the Leray- Schauder degree as one moves along the parameter in the continuity scheme. Let us observe that formula (8) appeared first in [10].

In contrast with the case $n = 5$ where there is only one blow up point, see Theorem 1.6, the cases $n \geq 6$ can present multiple blow up points, just as for scalar curvature problem in dimension $n \geq 4$, see [9] [11], [30]. In order to describe our results for $n = 6$, we introduce some more notation.
Given $f \in C^2(\mathbb{S}^6)$, consider
$$\mathcal{F} = \{p \in \mathbb{S}^6 \ : \ \nabla f(p) = 0\},$$
$$\mathcal{F}^+ = \{p \in \mathbb{S}^6 \ : \ \nabla f(p) = 0, \Delta_h f(p) > 0\},$$
$$\mathcal{F}^- = \{p \in \mathbb{S}^6 \ : \ \nabla f(p) = 0, \Delta_h f(p) < 0\},$$
$$\mathcal{M}_f = \{v \in C^2(\mathbb{S}^6) \ : \ v \text{ satisfies } (P)\}.$$

For $p_0 \in \mathbb{S}^6$, let $\pi_{p_0} : \mathbb{S}^6 \to \mathbb{R}^6$ denote the stereographic projection with pole $-p_0$. In $\pi_{p_0}$-stereographic coordinates, we consider the function $J_{p_0} : \mathbb{S}^6 \to \mathbb{R}$ given by

$$(9) \qquad J_{p_0}(y) = \frac{1}{2}\left(\frac{1+|y|^2}{|y|^2}\right), \quad y \in \mathbb{R}^6.$$

We note that the function $J_{p_0}$ is singular at $p_0$ and is a multiple of the Green's function of $P_h^6$ on $\mathbb{S}^6$. We recall that the Green's function for this operator with pole $p_0$ is a positive function $G_{p_0}$ satisfying $P_h^6 G_{p_0} = \delta_{p_0}$. One can check the existence of such a function taking a multiple of $J_{p_0}$, the uniqueness following from the coercivity of $P_h^6$, see above.
To each $\{p^1, \ldots, p^k\} \subseteq \mathcal{F} \setminus \mathcal{F}^-$, $k \geq 1$, we associate a $k \times k$ symmetric matrix $M = M(p^1, \ldots, p^k)$ defined by

$$(10) \qquad \begin{cases} M_{ii} = \frac{\Delta_h f(p^i)}{(f(p^i))^{\frac{3}{2}}}, \\ M_{ij} = -15 \frac{J_{p^i}(p^j)}{(f(p^i) f(p^j))^{\frac{1}{4}}} & \text{if } i \neq j. \end{cases}$$

Let $\rho = \rho(p^1, \ldots, p^k)$ denote the least eigenvalue of $M$. In particular, when $k = 1$ it is

$$\rho(p^1) = \frac{\Delta_h f(p^1)}{(f(p^1))^{\frac{3}{2}}}.$$

It has been first pointed out by A. Bahri, [8], see also [10], that when the interaction between different bubbles is of the same order as the 'self interaction', the function $\rho$ for a matrix similar to that given in (10) plays a fundamental role in the theory of the critical points at infinity. For Paneitz operator such kind of phenomenon appears when $n = 6$.
Define the set $\mathcal{A}$ to be

$$\mathcal{A} = \{f \in C^2(\mathbb{S}^6) \ : \ f \text{ is a positive Morse function on } \mathbb{S}^6 \text{ such that }$$
$$\Delta_h f \neq 0 \text{ on } \mathcal{F}, \text{ and } \rho = \rho(p^1, \ldots, p^k) \neq 0, \forall \, p^1, \ldots, p^k \in \mathcal{F}^+\}.$$

Let us observe that $\mathcal{A}$ is open in $C^2(\mathbb{S}^6)$ and dense in the space of positive $C^2$ function with respect to the $C^2$-norm. We introduce an integer valued continuous function Index : $\mathcal{A} \to \mathbb{N}$ by the following formula

$$\text{Index}(f) = -1 + \sum_{k=1}^{l} \sum_{\substack{\rho(p^{i_1}, \ldots p^{i_k}) > 0, \\ 1 \leq i_1 < i_2 \leq \cdots \leq i_k \leq l}} (-1)^{k-1+\sum_{j=1}^{k} m(f, p^{i_j})},$$

where $m(f, p^{i_j})$ denotes the Morse index of $f$ at $p^{i_j}$, and $l = card\,|\mathcal{F}^+|$.
Now we state our existence and compactness result for $n = 6$.

**Theorem 1.9** *Let $n = 6$ and assume that $f \in \mathcal{A}$. Then for all $\alpha \in (0, 1)$, there exists some constant $C$ depending only on $\min_{\mathbb{S}^6} f$, $\|f\|_{C^2(\mathbb{S}^6)}$, $\min_{\mathcal{F}} |\Delta_h f|$, $\min\{\rho(p^1, \ldots, p^k) : p^1, \ldots, p^k \in \mathcal{F}^+, k \geq 2\}$, and there exists some constant $\tilde{C} = \tilde{C}(C, \alpha)$ such that*

$$\frac{1}{C} \leq v \leq C, \quad \|v\|_{C^{4,\alpha}(\mathbb{S}^6)} \leq \tilde{C}, \tag{11}$$

*for all solution $v$ of equation $(P)$. Furthermore, for all $R \geq \max(C, \tilde{C})$*

$$deg(v - (P_h^n)^{-1}(f\,v^5), \mathcal{O}_R, 0) = Index(f), \tag{12}$$

*where*

$$\mathcal{O}_R = \left\{ v \in C^{4,\alpha}(\mathbb{S}^6) : \frac{1}{R} \leq v \leq R, \|v\|_{C^{4,\alpha}(\mathbb{S}^6)} \leq R \right\},$$

*and deg denotes the Leray-Schauder degree in $C^{4,\alpha}(\mathbb{S}^6)$. As a consequence, if $\text{Index}(f) \neq 0$, then $(P)$ has at least one solution.*

Theorem 1.9 is proved by using the following subcritical approximation of $(P)$

$$P_h^n u = \frac{n-4}{2} f\,u^{5-\tau}, \quad u > 0, \qquad \text{on } \mathbb{S}^6, \tag{13}$$

for $\tau > 0$ small. Using the Implicit Function Theorem, similarly to [38] and [30], we find for any $k$-tuple of points $p^1, \ldots, p^k \in \mathcal{F}^+$ with $\rho(p^1, \ldots, p^k) > 0$ a family of solutions of (13) highly concentrated around these points. Differently from the scalar curvature case, it is not immediate to check that this kind of solutions are positive: this fact is proved in Subsection 3.4. Using the blow up analysis we prove that the only blowing up solutions of (13) are of the above type. Then by a degree argument, using the condition $\text{Index}(f) \neq 0$, we derive the existence of at least one solution of $(P)$.

We organize our paper as follows. In Section 2 we provide the main local blow up analysis for (3). In section 3 we prove a-priori estimates in $H_2^2$ norm and $L^\infty$ norm for solutions of equation $(P)$ in dimensions 5 and 6. Then we derive our compactness and existence results. In the Appendix, we provide some useful technical results.

### Acknowledgements

Part of this work was accomplished when the first author enjoyed the hospitality of the SISSA at Trieste. He would like to thank A.Ambrosetti for the invitation and all the members of this institute for their hospitality. It is also a great pleasure to thank A.Ambrosetti, A.Bahri, A.Chang, M.Gursky ,E. Hebey and P.Yang for their interest in this work and their useful comments. The authors are also indebted to M.Morini for helpful discussions. A.M. is supported by a Fulbright Fellowship for the academic year 2000-2001 and by M.U.R.S.T. under the national project *"Variational methods and nonlinear differential equations"*. M.O.A. research is supported by a postdoctoral fellowship from SISSA.

## 2  Estimates for isolated simple blow up points

In this section we study the properties of isolated simple blow up points for equation (3). We first prove some Harnack type inequalities. In the following, given $r > 0$, $B_r$ will denote the open ball of radius $r$ centred at 0 in $\mathbb{R}^n$, and $\overline{B}_r$ its closure.

**Lemma 2.1** *For $\sigma \in (0, 1]$, let $A_\sigma = B_{3\sigma} \setminus \overline{B}_{\frac{1}{3}\sigma}$, and $A'_\sigma = \overline{B}_{2\sigma} \setminus B_{\frac{1}{2}\sigma}$. Let $g \in L^\infty(A_\sigma)$, and suppose $u$ is a positive $C^4$ solution of*

$$\Delta^2 u = g\, u, \quad \Delta u > 0, \quad \text{in } B_{3\sigma}.$$

*Then there exists $C = C(n, \|g\|_{L^\infty(A_\sigma)})$ such that*

$$\max_{A'_\sigma} u \leq C\, \overline{u}(2\sigma).$$

PROOF  Set $\xi(y) = u(\sigma y)$, $y \in B_3$. It is easy to see that $\xi$ satisfies

(14) $$\Delta^2 \xi(y) = \sigma^4\, g(\sigma y)\, \xi(y), \quad \xi(y) > 0, \quad \Delta \xi(y) > 0 \quad y \in B_3.$$

Set $w_1 = \xi$, $w_2 = \Delta \xi$. Then $(w_1, w_2)$ is a positive solution of the following elliptic system

(15) $$\begin{cases} \Delta w_1 = w_2 & \text{in } B_3; \\ \Delta w_2 = \sigma^4\, g(\sigma y)\, w_1 & \text{in } B_3. \end{cases}$$

For system (15), being a cooperative elliptic one, we can use the following weak Harnack Inequality due to Arapostathis, Ghosh and Marcus [3]

(16) $$\begin{cases} \max_{y \in \overline{B}_2 \setminus B_{\frac{1}{2}}} w_1(y) \leq C \max\left\{\min_{y \in \overline{B}_2 \setminus B_{\frac{1}{2}}} w_1(y), \min_{y \in \overline{B}_2 \setminus B_{\frac{1}{2}}} w_2(y)\right\} \\ \max_{y \in \overline{B}_2 \setminus B_{\frac{1}{2}}} w_2(y) \leq C \max\left\{\min_{y \in \overline{B}_2 \setminus B_{\frac{1}{2}}} w_1(y), \min_{y \in \overline{B}_2 \setminus B_{\frac{1}{2}}} w_2(y)\right\}, \end{cases}$$

where $C = C(n, \|g\|_{L^\infty(A_\sigma)})$. On the other hand, averaging the first equation in (15), we have

$$-\overline{\Delta w}_1(t) + \overline{w}_2(t) = 0, \quad \forall t \in (0, 3],$$

which is equivalent to

$$\ddot{\overline{w}}_1(t) + \frac{n-1}{t} \dot{\overline{w}}_1(t) + \overline{w}_2(t) = 0, \quad \forall t \in (0, 3].$$

Clearly, by the positivity of $w_1$ and $w_2$ and by (15), the functions $w_1$ and $w_2$ are superharmonic, so $\overline{w}_1$ and $\overline{w}_2$ are both non-negative and non-increasing. From another part, since $\xi$ is a smooth function, $\dot{\overline{w}}_1$ and $\dot{\overline{w}}_2$ are bounded near zero. Hence we can use the following Lemma

**Lemma 2.2** *(Serrin-Zou, [39]) Let $a$ be a positive real number and assume that $y = y(r) > 0$ satisfies*

$$\ddot{y}(t) + \frac{n-1}{t} \dot{y}(t) + \varphi(t) \leq 0 \quad \forall t \in (0, a),$$

*with $\varphi$ non-negative and non-increasing on $(0, a)$, and $\dot{y}$ bounded near $0$. Then there holds*

$$y(t) \geq c\, t^2\, \varphi(t) \quad \forall t \in (0, a),$$

*where $c = c(n)$.*

As a consequence, here, we have

$$\overline{w}_1(t) \geq c\, t^2\, \overline{w}_2(t) \quad \forall t \in (0, 2].$$

which implies that for $t \in [\frac{1}{2}, 2]$ it is $\overline{w}_2(t) \leq c_1\, \overline{w}_1(t)$, where $c_1$ is a positive real constant independent of $t$. Hence, using (16) and the superharmonicity of $w_1$, we deduce

$$\max_{\overline{B}_2 \setminus B_{\frac{1}{2}}} w_1 \leq C \max\left\{\min_{r \in [1/2, 2]} \overline{w}_1(r), \min_{r \in [1/2, 2]} \overline{w}_2(r)\right\} \leq C \min_{r \in [1/2, 2]} \overline{w}_1(r) \leq C\, \overline{w}_1(2).$$

This concludes the proof of Lemma 2.1, coming back to the original function $u$. ■

From Lemma 2.1, we deduce the following Harnack Inequality.

**Lemma 2.3** Let $\Omega \subseteq \mathbb{R}^n$ and $g \in L^\infty(\Omega)$, and assume that $u$ is a $C^4$ positive solution of

$$\Delta^2 u = g\,u, \quad \Delta u > 0, \quad \text{on } \Omega.$$

Then for every $\Omega' \subset\subset \Omega$, there exists $C = C(n, \|g\|_{L^\infty(\Omega)}, \Omega, \Omega')$ such that

$$\max_{\overline{\Omega'}} u \leq C \min_{\overline{\Omega'}} u.$$

PROOF  Let $\overline{\sigma} \in (0, \frac{1}{10} dist(\Omega', \partial\Omega))$. First we claim that there exists a constant $\tilde{C} = \tilde{C}(n, \|g\|_{L^\infty(\Omega)})$ such that for every point $x_0 \in \Omega'$ there holds

(17) $$\max_{\overline{B}_\sigma(x_0)} u \leq \tilde{C}\, u(x_0).$$

Once (17) is established, the assertion follows by covering $\Omega'$ with a finite number of balls of radius $\overline{\sigma}$ starting from a minimum of $u$ on $\Omega'$. Of course, the number of these balls can be chosen depending only on $n, \Omega, \Omega'$.

Let us now prove (17). Consider the function $z(y) = u(x_0 + y)$. Then it is clear that $z$ satisfies the assumptions of Lemma 2.1 for $\sigma \in (0, \overline{\sigma})$, and taking $g(x_0 + \cdot)$ instead of $g$. Hence we deduce

$$\max_{\overline{\Omega'}} z \leq \tilde{C}(n, \|g\|_{L^\infty(\Omega)})\, \overline{z}(2\sigma), \quad \sigma \in (0, \overline{\sigma}).$$

Recalling the definition of $z$, and taking into account that $z$ is superharmonic, we have

$$\max_{\overline{\Omega'}} z \leq \tilde{C}(n, \|g\|_{L^\infty(\Omega)})\, \overline{z}(0) = \tilde{C}(n, \|g\|_{L^\infty(\Omega)})\, u(x_0).$$

This implies (17). ∎

**Remark 2.4** Let $f$ be a $C^1$ positive function on $\mathbb{S}^n$, and $q \in [1, \frac{n+4}{n-4}]$. Let $u$ be a positive solution of $P_h^n u = f\,u^q$ on $\mathbb{S}^n$. It follows from Lemma 2.3 that upper bounds on $u$ imply also lower bounds on $u$.

**Lemma 2.5** Let $\{f_i\}_i$ satisfy (2), $\{u_i\}_i$ satisfy (3), and let $y_i \to \overline{y} \in \Omega$ be an isolated blow up point. Then for any $r \in (0, \frac{1}{3}\overline{r})$, we have the following weak Harnack inequality

$$\max_{y \in \overline{B}_{2r}(y_i) \setminus B_{\frac{r}{2}}(y_i)} u_i(y) \leq C \min_{y \in \overline{B}_{2r}(y_i) \setminus B_{\frac{r}{2}}(y_i)} u_i(y)$$

$$\max_{y \in \overline{B}_{2r}(y_i) \setminus B_{\frac{r}{2}}(y_i)} (\Delta u_i(y)) \leq \frac{1}{r^2} C \min_{y \in \overline{B}_{2r}(y_i) \setminus B_{\frac{r}{2}}(y_i)} u_i(y),$$

where $C$ is some positive constant depending only on $n$, $\overline{C}$, and $\sup_i \|f_i\|_{L^\infty(B_{\overline{r}}(y_i))}$.

PROOF  Set $\xi(y) = r^{\frac{4}{q_i - 1}} u_i(y_i + r\,y)$, $y \in B_3$. It is easy to see that $\xi$ satisfies

(18) $$\begin{cases} \Delta^2 \xi(y) = \frac{n-4}{2} f_i(y_i + r\,y)\xi(y)^{q_i} & y \in B_3, \\ \Delta \xi > 0 & y \in B_3, \\ 0 < \xi(y) \leq \overline{C}\,|y|^{-\frac{4}{q_i - 1}} & y \in B_3. \end{cases}$$

The first inequality follows easily from Lemma 2.3. For the second one, it is sufficient to use the above rescaling, (16) and Lemma 2.2. ∎

**Remark 2.6** *It is clear from the proof of Lemma 2.5 that the conclusion remains true if instead of assuming that $\overline{y}$ is an isolated blow up we only assume* (4).

**Proposition 2.7** *Let $\{f_i\}_i$ be bounded in $C^1_{loc}(\Omega)$ and satisfy* (2). *Let $\{u_i\}_i$ satisfy* (3), *and let $y_i \to \overline{y} \in \Omega$ be an isolated blow up point for $\{u_i\}_i$. Then, for any $R_i \to +\infty$ and $\varepsilon_i \to 0^+$, we have, after passing to a subsequence of $u_i$ (still denoted $\{u_i\}_i$, $\{y_i\}_i$, etc...), that*

$$\left\| u_i(y_i)^{-1} u_i \left( u_i(y_i)^{-\frac{q_i-1}{4}} \cdot + y_i \right) - (1 + k_i |\cdot|^2)^{\frac{4-n}{2}} \right\|_{C^4(B_{2R_i}(0))} \leq \varepsilon_i,$$

$$R_i \, u_i(y_i)^{-\frac{q_i-1}{4}} \to 0 \quad \text{as } i \to +\infty,$$

*where $k_i^2 = \frac{1}{2n\,(n-2)\,(n+2)} f_i(y_i)$.*

PROOF Consider

$$\xi_i(x) = u_i(y_i)^{-1} u_i \left( u_i(y_i)^{-\frac{q_i-1}{4}} x + y_i \right) \quad \forall |x| \leq \overline{r} \, u_i(y_i)^{\frac{q_i-1}{4}}.$$

Clearly

$$\begin{cases} \Delta^2 \xi_i(x) = \frac{n-4}{2} f_i \left( u_i(y_i)^{-\frac{q_i-1}{4}} x + y_i \right) \xi_i(x)^{q_i} & |x| \leq \overline{r} \, u_i(y_i)^{\frac{q_i-1}{4}}, \\ 0 < \xi_i(x) \leq \overline{C} \, |x|^{-\frac{4}{q_i-1}}, \quad \Delta \xi_i > 0 & |x| \leq \overline{r} \, u_i(y_i)^{\frac{q_i-1}{4}}, \\ \xi_i(0) = 1, \quad \nabla \xi_i(0) = 0. \end{cases}$$

It follows from Remark 2.6 and from the superharmonicity of $\xi_i$ that for $r \in (0, +\infty)$ we have for $i$ large

$$\max_{x \in \partial B_r(0)} \xi_i(x) \leq C \min_{x \in \partial B_r(0)} \xi_i(x) \leq C \, \xi_i(0) = C,$$

so $\xi_i$ is uniformly bounded in $C^\infty_{loc}(\mathbb{R}^n)$. For every $r > 1$, by Remark 2.6, we also have

(19) $$\max_{\partial B_r} \Delta \xi \leq \frac{C}{r^2} \overline{\xi}_i(r) \leq C.$$

Since the functions $\Delta \xi_i$ satisfy the equation $\Delta(\Delta \xi_i) = f_i(y_i + r\, y) \xi_i^{q_i}$, then from $L^p$ estimates (see e.g. [28], Theorem 9.11) and Schauder estimates (see e.g. [28], chapter 6) $\{\Delta \xi_i\}_i$ is bounded in $C^{2,\alpha}_{loc}(\mathbb{R}^n)$. By the same reasons it follows that $\{\xi_i\}_i$ is bounded in $C^{4,\alpha}_{loc}(\mathbb{R}^n)$.

Hence by the Ascoli-Arzelà Theorem, there exists some function $\xi$ such that, after passing to a subsequence,

$$\begin{cases} \xi_i \to \xi & \text{in } C^4_{loc}(\mathbb{R}^n), \\ \Delta^2 \xi = \frac{n-4}{2} (\lim_i f_i(y_i)) \, \xi^{\frac{n+4}{n-4}} & \text{in } \mathbb{R}^n, \\ \xi \geq 0, \quad \Delta \xi_i \geq 0, \quad \xi(0) = 1, \quad \nabla \xi(0) = 0. \end{cases}$$

Since $\Delta \xi \geq 0$, and $\xi \geq 0$, it follows from the maximum principle that $\xi$ is positive in $\mathbb{R}^n$. It follows from standard regularity arguments that $\xi$ is $C^\infty$ in $\mathbb{R}^n$, so the result in Lin [31] implies that

$$\xi(x) = \left( 1 + \lim_i k_i |x|^2 \right)^{\frac{4-n}{2}}.$$

*where $k_i^2 = \frac{1}{2n\,(n-2)\,(n+2)} f_i(y_i)$. Proposition 2.7 is now proved.* ∎

**Remark 2.8** *It follows from the proof of Proposition 2.7 that, under the assumption that $\{f_i\}_i$ is bounded in $C^1_{loc}(\Omega)$, if a sequence of solutions $\{u_i\}$ of* (3) *remains bounded in $L^\infty_{loc}(\Omega)$, then $u_i$ converges in $C^{4,\alpha}_{loc}$ along a subsequence.*

Since passing to subsequences does not affect our proofs, we will always choose $R_i \to +\infty$ first, and then $\varepsilon_i$ (depending on $R_i$) as small as necessary. In particular $\varepsilon_i$ is chosen small enough so that $y_i$ is the only critical point of $u_i$ in $0 < |y| < R_i \, u_i(y_i)^{-\frac{q_i-1}{4}}$, $\hat{u}_i(r)$ (defined after formula (5)) has a unique critical point in $\left(0, R_i \, u_i(y_i)^{-\frac{q_i-1}{4}}\right)$, $2\varepsilon_i < \left(1+k_i R_i^2\right)^{\frac{4-n}{2}}$, and

$$\left\| u_i(y_i)^{-1} \, u_i \left( u_i(y_i)^{-\frac{q_i-1}{4}} \cdot + y_i \right) - (1 + k_i |\cdot|^2)^{\frac{4-n}{2}} \right\|_{C^4(B_{2R_i}(0))} \to 0.$$

**Proposition 2.9** *Let $\{f_i\}_i \in C^1_{loc}(B_2)$ satisfies (2) with $\Omega = B_2$ and*

(20) $$|\nabla f_i(y)| \leq A_2, \quad \forall y \in B_2,$$

*for some positive constant $A_2$. Assume that $\{u_i\}_i$ satisfies (3) with $\Omega = B_2$, and let $y_i \to \overline{y} \in \Omega$ be an isolated simple blow up point for $\{u_i\}_i$, which for some positive constant $A_3$ satisfies*

(21) $$|y - y_i|^{\frac{4}{q_i-1}} u_i(y_i) \leq A_3, \quad \forall y \in B_2.$$

*Then there exists some positive constant $C = C(n, A_1, A_2, A_3, \varrho)$ ($\varrho$ being given in the definition of isolated simple blow up point) such that for $R_i \, u_i(y_i)^{-\frac{q_i-1}{4}} \leq |y - y_i| \leq 1$*

(22) $$C^{-1} u_i(y_i)^{-1} |y-y_i|^{4-n} \leq u_i(y) \leq C \, u_i(y_i)^{-1} |y-y_i|^{4-n}.$$

*Furthermore there exists some biharmonic function $b(y)$ in $B_1$ such that we have, after passing to a subsequence,*

$$u_i(y_i) \, u_i(y) \to h(y) = a \, |y|^{4-n} + b(y) \quad \text{in } C^4_{loc}(B_1 \setminus \{0\}),$$

*where*

$$a = \left(\lim_i k_i\right)^{\frac{4-n}{2}}.$$

Before proving Proposition 2.9 we need some Lemmas.

**Lemma 2.10** *Under the assumptions of Proposition 2.9, except for (20), there exist $\delta_i > 0$, $\delta_i = O(R_i^{-4+o(1)})$ such that*

$$u_i(y) \leq C \, u_i(y_i)^{-\lambda_i} |y - y_i|^{4-n+\delta_i} \quad \text{for } R_i \, u_i(y_i)^{-\frac{q_i-1}{4}} \leq |y-y_i| \leq 1,$$

*where $\lambda_i = (n - 4 - \delta_i) \left(\frac{q_i-1}{4}\right) - 1$.*

PROOF Let $r_i = R_i \, u_i(y_i)^{-\frac{q_i-1}{4}}$; it follows from Proposition 2.7 that

(23) $$u_i(y) \leq C \, u_i(y_i) \, R_i^{4-n}, \quad \Delta u_i(y) \leq C \, u_i(y_i)^{\frac{q_i+1}{2}} R_i^{2-n} \quad \text{for all } |y - y_i| = r_i.$$

Let $\hat{u}_i(r)$ be defined as in (5); it follows from the definition of isolated simple blow up and Proposition 2.7 that there exists $\varrho > 0$ such that

(24) $$r^{\frac{4}{q_i-1}} \overline{u}_i(r) \text{ is strictly decreasing for } r_i < r < \varrho.$$

Using (23), (24) and Lemma 2.5 we have that for all $r_i < |y - y_i| < \varrho$

$$|y-y_i|^{\frac{4}{q_i-1}} u_i(y) \leq C \, |y-y_i|^{\frac{4}{q_i-1}} \overline{u}_i(|y-y_i|) \leq C \, r_i^{\frac{4}{q_i-1}} \overline{u}_i(r_i) \leq C \, R_i^{\frac{4-n}{2}+o(1)}.$$

Therefore
$$u_i(y)^{q_i-1} = O\left(R_i^{-4+o(1)} |y-y_i|^{-4}\right) \quad \text{for all } r_i \le |y-y_i| \le \varrho. \tag{25}$$

Now we would like to apply Lemma 4.3 with $D = \{r_i \le |y-y_i| \le \varrho\}$, and $L_1 = L_2 = -\Delta$, $h_{11} = 0$, $h_{12} = 1$, $h_{21} = \frac{n-4}{2} f_i u_i^{q_i-1}$.

Take $\alpha \in \left(0, \frac{1}{n-4}\right)$, and let
$$\varphi_1 = |y-y_i|^{-\alpha}, \quad \varphi_2 = \Delta \varphi_1.$$

By a direct computation, taking into account (25), one can check that
$$\begin{cases} -\Delta \varphi_1 + \varphi_2 = 0; \\ -\Delta \varphi_2 + \frac{n-4}{2} f_i u_i^{q_i-1} = \left[-\alpha(2+\alpha)(n-2-\alpha)(n-4-\alpha) + O\left(R_i^{-4+o(1)}\right)\right] |y-y_i|^{-(4+\alpha)}, \end{cases}$$

for $r_i \le |y-y_i| \le \varrho$. It can be easily seen that with our choice of $\alpha$ it is $-\Delta \varphi_2 + \frac{n-4}{2} f_i u_i^{q_i-1} < 0$.
Now set
$$\varphi_1 = |y-y_i|^{4-n+\delta_i}, \quad \varphi_2 = \Delta \varphi_1.$$

Then there holds, again by (25)
$$\begin{cases} -\Delta \varphi_1 + \varphi_2 = 0, \\ -\Delta \varphi_2 + \frac{n-4}{2} f_i u_i^{q_i-1} = \left[-\delta_i(2+\delta_i)(n-2-\delta_i)(n-4-\delta_i) + O\left(R_i^{-4+o(1)}\right)\right] |y-y_i|^{-n+\delta_i}, \end{cases}$$

for $r_i \le |y-y_i| \le \varrho$. So we can choose $\delta_i = O(R_i^{-4+o(1)})$ such that $-\Delta \varphi_2 + \frac{n-4}{2} f_i u_i^{q_i-1} < 0$. Now set
$$\psi_i(y) = \gamma_1 M_i \varrho^\alpha |y-y_i|^{-\alpha} + \gamma_2 u_i(y_i)^{-\lambda_i} |y-y_i|^{4-n+\delta_i},$$

where $M_i = \max_{\partial B_\varrho} u_i$, $\lambda_i = (n-4-\delta_i)\left(\frac{q_i-1}{4}\right) - 1$ and $\gamma_1, \gamma_2 > 0$. It follows from the previous computations that we can apply Lemma 4.3 with $(w_1, w_2) = (\psi_i, \Delta \psi_i)$ and $(z_1, z_2) = (u_i, \Delta u_i)$ provided we show
$$\begin{cases} u_i \le \psi_i & \text{on } \partial(\{r_i \le |y-y_i| \le \varrho\}), \\ \Delta u_i \le \Delta \psi_i & \text{on } \partial(\{r_i \le |y-y_i| \le \varrho\}). \end{cases} \tag{26}$$

For this purpose we observe that for $|y-y_i| = \varrho$, it is $\psi_i \ge \gamma_1 M_i$, so if $\gamma_1 > 1$ we have
$$\psi_i(y) \ge u_i(y) \quad \text{for } |y-y_i| = \varrho.$$

Moreover, by Lemma 2.5 there exist $C > 0$ such that
$$\max_{\partial B_\varrho} \Delta u_i \le C \overline{u}_i(\varrho) \le C M_i,$$

so one can easily check that for some $\gamma_1 > 0$ sufficiently large there holds
$$\Delta \psi_i(y) \ge \Delta u_i(y) \quad \text{for } |y-y_i| = \varrho.$$

We observe that we have proved (26) on $|y-y_i| = \varrho$; for $|y-y_i| = r_i$, we have
$$\psi_i(y) \ge \gamma_2 u_i(y_i)^{-\lambda_i} r_i^{4-n+\delta_i}.$$

But $r_i = R_i u_i(y_i)^{-\frac{q_i-1}{4}}$ so, taking into account the expression of $\lambda_i$ we derive
$$\psi_i(y) \ge \gamma_2 u_i(y_i) R_i^{4-n+\delta_i} \quad \text{for } |y-y_i| = r_i.$$

By Proposition 2.7, it turns out that $u_i(y) \leq C\, u_i(y_i)\, R_i^{4-n}$ for $|y - y_i| = r_i$, so it follows that for $i$ large
$$\psi_i(y) \geq u_i(y) \quad \text{for } |y - y_i| = r_i.$$

From another part, it is
$$\Delta \psi_i(y) \geq \gamma_2\, u_i(y_i)^{-\lambda_i}\, (n - 4 + \delta_i)\, (2 + \delta_i)\, r_i^{2-n+\delta_i} \quad \text{for } |y - y_i| = r_i,$$

so from the expression of $r_i$ and from (23) it follows that for $\gamma_2$ large enough
$$\Delta \psi_i(y) \geq \Delta u_i(y) \quad \text{for } |y - y_i| = r_i.$$

We have now proved (26), so it is:

(27) $$u_i(y) \leq \psi_i(y) \quad \text{for } r_i \leq |y - y_i| \leq \varrho.$$

Now, since $y_i \to \overline{y}$ is an isolated simple blow up, taking into account (24), Lemma 2.5 and inequality (27), we deduce that for $r_i < \theta < \varrho$ it is

$$\begin{aligned}
\varrho^{\frac{q_i-1}{4}} M_i &\leq C\, \varrho^{\frac{q_i-1}{4}} \overline{u}_i(\varrho) \leq C \theta^{\frac{q_i-1}{4}} \overline{u}_i(\theta) \\
&\leq C \theta^{\frac{q_i-1}{4}} \left( \gamma_1 M_i\, \varrho^\alpha\, \theta^{-\alpha} + \gamma_2\, u_i(y_i)^{-\lambda_i} \theta^{4-n+\delta_i} \right).
\end{aligned}$$

Since we are assuming $0 < \alpha < \frac{2}{n-4}$, we can choose $\theta$ independent of $i$ such that
$$C \gamma_1\, \varrho^\alpha\, \theta^{\frac{q_i-1}{4} - \alpha} < \frac{1}{2}\, \varrho^{\frac{q_i-1}{4}},$$

and with such a choice it turns out that

(28) $$M_i \leq C\, u_i(y_i)^{-\lambda_i}.$$

This concludes the proof of the Lemma for $r_i \leq |y - y_i| \leq \varrho$; for $\varrho \leq |y - y_i| \leq 1$, it is sufficient to use Lemma 2.5. ∎

We recall that we have set $\tau_i = \frac{n+4}{n-4} - q_i$.

**Lemma 2.11** *Under the hypotheses of Proposition 2.9 we have*
$$\tau_i = O\left( u_i(y_i)^{-\frac{2}{n-4}+o(1)} \right),$$

*and therefore*
$$u_i(y_i)^{\tau_i} = 1 + o(1).$$

PROOF  Applying Proposition 4.1 with $r = 1$ we obtain
$$\frac{n-4}{2(q+1)} \sum_{j=1}^n \int_{B_1} x_j \frac{\partial f_i}{\partial x_j} u_i^{q_i+1} dx + \frac{n-4}{2} \left( \frac{n}{q_i+1} - \frac{n-4}{2} \right) \int_{B_1} f_i\, u_i^{q_i+1}\, dx$$
$$- \frac{n-4}{2(q_i+1)} \int_{\partial B_1} f_i\, u_i^{q_i+1} d\sigma = \int_{\partial B_1} B(r, x, u, \nabla u_i, \nabla^2 u_i, \nabla^3 u_i)\, d\sigma.$$

From (28), Lemma 2.5 and from standard elliptic estimates, one can easily deduce that
$$\int_{\partial B_1} B(r, x, u_i, \nabla u, \nabla^2 u_i, \nabla^3 u_i)\, d\sigma = O\left( u_i(y_i)^{-2+o(1)} \right), \quad \int_{\partial B_1} f_i\, u_i^{q_i+1} d\sigma = O\left( u_i(y_i)^{-\frac{2n}{n-4}+o(1)} \right).$$

Moreover, using Proposition 2.7 and simple rescaling arguments we derive

$$\sum_{j=1}^{j=n} \int_{B_1} x_j \frac{\partial f_i}{\partial x_j} u_i^{q_i+1} \, dx = O\left(u_i(y_i)^{-\frac{2}{n-4}+o(1)}\right).$$

Hence it follows that

$$\tau_i = O\left(u_i(y_i)^{-\frac{2}{n-4}+o(1)}\right) + O\left(u_i(y_i)^{-2+o(1)}\right) + O\left(u_i(y_i)^{-\frac{2n}{n-4}+o(1)}\right) = O\left(u_i(y_i)^{-\frac{2}{n-4}+o(1)}\right).$$

This concludes the proof. ∎

**Lemma 2.12** *There holds*

$$u_i(y_i) \, u_i^{q_i}(y) \to l \, \delta_0(y),$$

*where* $l = \frac{2}{n(n+2)} \left(\lim_i k_i\right)^{-\frac{n}{2}} \omega_{n-1}$ *and the convergence is in the weak sense of measures.*

PROOF   Take $\psi \in C_c^\infty(B_1)$; we clearly have

$$\int_{B_1} u_i(y_i) \, u_i^{q_i} \, \psi = \int_{B_{r_i}} u_i(y_i) \, u_i^{q_i} \, \psi + \int_{B_1 \setminus B_{r_i}} u_i(y_i) \, u_i^{q_i} \, \psi.$$

Using Proposition 2.7 we deduce by simple computations

(29) $$\int_{B_{r_i}} u_i(y_i) \, u_i^{q_i} \, \psi = (\psi(0) + o(1)) \, u_i(y_i) \int_{B_{r_i}} u_i^{q_i} \to l \, \psi(0).$$

Moreover, by Lemma 2.10 there holds

$$\int_{B_1 \setminus B_{r_i}} u_i^{q_i} \leq C \int_{B_1 \setminus B_{r_i}} \left(u_i(y_i)^{-\lambda_i} |y - y_i|^{4-n+\delta_i}\right)^{q_i}$$

(30) $$\leq C R_i^{n-q_i(n-4-\delta_i)} \, u_i(y_i)^{-1+O(\tau_i)} = o(1) \, u_i(y_i)^{-1},$$

so the conclusion follows. ∎

**Lemma 2.13** *Let* $w_i : B_1 \to \mathbb{R}$ *be defined as*

$$w_i(y) = u_i(y_i) \, \Delta \, u_i(y),$$

*and let* $H_{B_1}(w_i)$ *denote the unique function satisfying*

$$\begin{cases} \Delta \, H_{B_1}(w_i) = 0 & \text{in } B_1, \\ H_{B_1}(w_i) \equiv w_i & \text{on } \partial B_1. \end{cases}$$

*If we set* $\tilde{w}_i = w_i - H_{B_1}(w_i)$, *then we have*

$$\tilde{w}_i \to \tilde{l} \, G_{B_1}(0, x) \quad \text{in } L^1(B_1),$$

*where*

$$\tilde{l} = \frac{n-4}{n(n+2)} \left(\lim_i f_i(y_i)\right) \left(\lim_i k_i\right)^{-\frac{n}{2}} \omega_{n-1},$$

*and where* $G_{B_1}$ *denotes the Green's function of* $\Delta$ *in* $B_1$ *under Dirichlet boundary conditions.*

PROOF  It is easy to prove, using Lemma 2.12 and the Green's representation formula, that $\tilde{w}_i(\cdot) \to l\, G_{B_1}(0,\cdot)$ pointwise. We will prove that $\tilde{w}_i$ is bounded in $W_0^{1,q}(B_1)$ for $q < \frac{n}{n-1}$. Then the Lemma will follow from the Rellich compactness Theorem.

Hence we are reduced to prove that $\nabla \tilde{w}_i$ is bounded in $L^q(B_1)$ for all $q < \frac{n}{n-1}$. Let $f \in (L^{q'}(B_1))^n$, $q' > n$, and consider the problem

$$\begin{cases} -\Delta v = \mathrm{div}\, f & \text{in } W^{-1,q'}(B_1), \\ v \in W_0^{1,q'}(B_1). \end{cases}$$

This problem admits a unique solution $v$, which satisfies

$$\|v\|_{W_0^{1,q'}(B_1)} \leq C\, \|\mathrm{div} f\|_{W_0^{-1,q'}(B_1)} \leq C'\, \|f\|_{L^{q'}(B_1)}.$$

Moreover, since $q' > n$, $v \in C^0(B_1)$ and

(31) $$\|v\|_{C^0(B_1)} \leq C \|v\|_{W_0^{1,q'}(B_1)} \leq C'\, \|f\|_{L^{q'}(B_1)}.$$

We can split $f$ into

$$f = \nabla v + f_1 \quad \text{with } \mathrm{div}\, f_1 = 0;$$

hence we have

$$\left| \int_{B_1} f \cdot \nabla \tilde{w}_i \right| = \left| \int_{B_1} \nabla v \cdot \nabla \tilde{w}_i + \int_{B_1} f_1 \cdot \nabla \tilde{w}_i \right| = \left| \int_{B_1} \nabla v \cdot \nabla \tilde{w}_i + \int_{B_1} \langle \mathrm{div} f_1, \tilde{w}_i \rangle \right| = \left| \int_{B_1} \nabla v \cdot \nabla \tilde{w}_i \right|.$$

Writing $\int_{B_1} \nabla v \cdot \nabla \tilde{w}_i$ as $\int_{B_1} \Delta \tilde{w}_i\, v$, we deduce, using (31) and Lemma 2.12

$$\left| \int_{B_1} f \cdot \nabla \tilde{w}_i \right| \leq C'' \|v\|_{C^0(B_1)} \leq C\, \|f\|_{L^{q'}(B_1)}.$$

Hence by Banach-Steinhaus Theorem the sequence $\{\tilde{w}_i\}_i$ is bounded in $W^{1,q}$ for all $q < \frac{n}{n-1}$. This concludes the proof. ∎

**Lemma 2.14** *Let $\Gamma : (0,1) \to \mathbb{R}$ be defined by*

$$\Gamma(\sigma) = \lim_i \left( \int_{B_\sigma} \Delta u_i \right) \left( \int_{B_1 \setminus B_\sigma} \Delta u_i \right)^{-1}.$$

*Then, after extracting a subsequence if necessary, $\Gamma(\cdot)$ is well defined and satisfies $\Gamma(t) > 0\ \forall t \in (0,1)$ and $\Gamma(t) \to 0$ as $t \to 0$.*

PROOF  Let $w_i = u_i(y_i)\, \Delta u_i$. Since $H_{B_1}(w_i)$ is harmonic, by the mean value property, there holds

$$\int_{B_\sigma} H_{B_1}(w_i) = \omega_{n-1} \overline{H}_i \int_0^\sigma r^{n-1}, \quad \int_{B_1 \setminus B_\sigma} H_{B_1}(w_i) = \omega_{n-1} \overline{H}_i \int_\sigma^1 r^{n-1},$$

where we have set

$$\overline{H}_i = \frac{1}{\omega_{n-1}} \int_{\partial B_1} H_{B_1}(w_i).$$

In particular recalling the definition of $\tilde{w}_i$, it follows that

(32) $$\int_{B_\sigma} w_i = \int_{B_\sigma} \tilde{w}_i + \omega_{n-1} \overline{H}_i \int_0^\sigma r^{n-1}; \quad \int_{B_1 \setminus B_\sigma} w_i = \int_{B_1 \setminus B_\sigma} \tilde{w}_i + \omega_{n-1} \overline{H}_i \int_\sigma^1 r^{n-1}.$$

Now two cases may occur: the first is $\sup_i \overline{H}_i < +\infty$, and the second is $\limsup_i \overline{H}_i = +\infty$. We treat the two cases separately.

**Case 1**  $\sup_i \overline{H}_i < +\infty$

Taking into account (32) we have

$$\left(\int_{B_\sigma} w_i\right) \left(\int_{B_1 \setminus B_\sigma} w_i\right)^{-1} = \left(\int_{B_\sigma} \tilde{w}_i + \omega_{n-1} \overline{H}_i \int_0^\sigma r^{n-1}\right) \left(\int_{B_1 \setminus B_\sigma} \tilde{w}_i + \omega_{n-1} \overline{H}_i \int_\sigma^1 r^{n-1}\right)^{-1}.$$

By Lemma 2.13 it is $\tilde{w}_i \to \tilde{w} = l\, G_{B_1}(0, \cdot) > 0$ in $L^1(B_1)$ so, since up to a subsequence $\overline{H}_i \to \overline{H}$, $\Gamma(\sigma)$ is

$$\Gamma(\sigma) = \left(\int_{B_\sigma} \tilde{w} + \omega_{n-1} \overline{H} \int_0^\sigma r^{n-1}\right) \left(\int_{B_1 \setminus B_\sigma} \tilde{w} + \omega_{n-1} \overline{H} \int_\sigma^1 r^{n-1}\right)^{-1}.$$

**Case 2**  $\limsup_i \overline{H}_i = +\infty$

Arguing as in the first case, taking into account the convergence of $\tilde{w}_i$ and the fact that $\overline{H}_i \to +\infty$, we derive

$$\lim_i \left(\int_{B_\sigma} w_i\right) \left(\int_{B_1 \setminus B_\sigma} w_i\right)^{-1} = \left(\int_0^\sigma r^{n-1}\right) \left(\int_\sigma^1 r^{n-1}\right)^{-1}.$$

So in this case $\Gamma(\sigma)$ is

$$\Gamma(\sigma) = \left(\int_0^\sigma r^{n-1}\right) \left(\int_\sigma^1 r^{n-1}\right)^{-1}.$$

In both the cases, the function $\Gamma(\sigma)$ satisfies the required properties, hence the conclusion follows. ∎

**Lemma 2.15** *Set $v_i(y) = \overline{u}_i(1)^{-1} u_i(y)$. Then there holds*

$$v_i \to v(x) = a_1 |x|^{4-n} + b(x) \quad \text{in } C^2_{loc}(B_1 \setminus \{0\}),$$

*where*

$$a_1 > 0, \quad v > 0, \quad b \in C^\infty(B_1), \quad \Delta^2 b = 0.$$

PROOF   It is easy to see that $v_i$ satisfies

$$\Delta^2 v_i = \frac{n-4}{2} \overline{u}_i(1)^{q_i-1} f_i(y) v_i^{q_i} \quad \text{in } B_2.$$

We observe that by Lemma 2.10, $\overline{u}_i(1) \to 0$ so it follows from Lemma 2.5 and standard elliptic estimates (we note that $\overline{v}_i(1) = 1$ as $i \to +\infty$) that $\{v_i\}$ converges in $C^4_{loc}(B_1 \setminus \{0\})$ to some function $v$ which satisfies

$$\Delta^2 v = 0, \quad v \geq 0, \quad \Delta v \geq 0 \quad \text{in } B_2 \setminus \{0\}.$$

Moreover $v$ must possess a singularity at 0. In fact, since we are assuming that $y_i \to \overline{y}$ is an isolated simple blow up, it follows that also $r^{\frac{n-4}{2}} \overline{v}(r)$ is non-increasing for $0 < r < \varrho$, which is impossible if $v$ is regular near the origin.

It follows from Lemma 4.5 that $v$ is of the form

$$v(x) = a_1 |x|^{4-n} + a_2 |x|^{2-n} + b(x),$$

where $a_1, a_2 \geq 0$, and $b \in C^\infty(B_1)$ satisfies $\Delta^2 b = 0$. Since $v$ is singular at 0, it is $a_1 + a_2 > 0$. Using the divergence Theorem and the convergence of $v_i$ to $v$, we derive that for $\sigma \in (0,1)$

$$\lim_i \overline{u}_i(1)^{-1} \int_{B_\sigma} \Delta u_i = \int_{\partial B_\sigma} -\frac{\partial v}{\partial \nu} = \omega_{n-1} a_2 (n-2) + o_\sigma(1),$$

where $o_\sigma(1) \to 0$ as $\sigma \to 0$. Using Lemma 2.14 we deduce

$$\omega_{n-1} a_2 (n-2) + o_\sigma(1) = \int_{\partial B_\sigma} -\frac{\partial v}{\partial \nu} = \frac{\int_{B_\sigma} \Delta u_i}{\int_{B_1 \setminus B_s} \Delta u_i} \lim_i \overline{u}_i(1)^{-1} \int_{B_1 \setminus B_\sigma} \Delta u_i \to \Gamma(\sigma) \int_{B_1 \setminus B_\sigma} \Delta v.$$

Since $\Delta v \in L^1(B_1)$, letting $\sigma \to 0$ we obtain $a_2 = 0$, taking into account that $\Gamma(\sigma) \to 0$ as $\sigma \to 0$. ∎

PROOF OF PROPOSITION 2.9. Let $v_i$ be as in Lemma 2.15. We prove first the inequality (22) for $|y - y_i| = 1$, namely we show that for some $C > 0$ independent of $i$

(33) $$C^{-1} \leq \overline{u}_i(1) \, u_i(y_i) \leq C.$$

Multiply (3) by $\overline{u}_i(1)^{-1}$ and integrate on $B_1$, we have

(34) $$\int_{\partial B_1} -\frac{\partial}{\partial \nu}(\Delta v_i) = \frac{n-4}{2} \overline{u}_i(1)^{-1} \int_{B_1} f_i u_i^{q_i},$$

where we have set, as before, $v_i(y) = \overline{u}_i(1)^{-1} u_i(y)$. Applying Lemma 2.15, we deduce using the biharmonicity of $b$

(35) $$\lim_i \int_{\partial B_1} -\frac{\partial}{\partial \nu}(\Delta v_i) = \int_{\partial B_1} \frac{\partial}{\partial \nu} \left(2 a_1 (4-n)|x|^{2-n} - \Delta b\right) = 2 a_1 (n-2)(n-4)\omega_{n-1} > 0.$$

From (34) and (35) we deduce that

(36) $$\lim_i \overline{u}_i(1)^{-1} \int_{B_1} f_i u_i^{q_i} = 4 a_1 (n-2) \omega_{n-1} > 0.$$

From another part, it follows from Lemma 2.12 that

(37) $$\lim_i u_i(y_i) \int_{B_1} f_i u_i^{q_i} \to l \left(\lim_i f_i(y_i)\right).$$

Hence (36) and (37) imply (33). To establish (22) for $r_i \leq |y - y_i| \leq 1$, we only need to scale the problem and to reduce it to the case $|y - y_i| = 1$. In fact, arguing by contradiction, if there exists a subsequence $\{\tilde{y}_i\}_i$ satisfying $r_i \leq |\tilde{y}_i - y_i| \leq \rho$, and $\lim_i u_i(\tilde{y}_i) u_i(y_i) |\tilde{y}_i - y_i|^{4-n} = \infty$, set $\tilde{r}_i = |\tilde{y}_i - y_i|$, $\tilde{u}_i(y) = (\tilde{r}_i)^{\frac{4}{q_i-1}} u_i(y_i + \tilde{r}_i y)$. Then $\tilde{u}_i$ satisfies all the assumptions of Proposition 2.9 with the same constants and with 0 instead of $\overline{y}$. It follows from (33) that

$$C^{-1} \leq \tilde{u}_i(0) \, \tilde{u}_i \left(\frac{\tilde{y}_i - y_i}{\tilde{r}_i}\right) \leq C.$$

This leads to a contradiction, so we have (22).

Next we compute the value of $a$. Multiplying equation (3) by $u_i(y_i)$ and integrating on $B_1$, we have

$$\int_{\partial B_1} -\frac{\partial}{\partial \nu}(u_i(y_i) \Delta u_i) = \frac{n-4}{2} u_i(y_i) \int_{B_1} f_i u_i^{q_i}.$$

It follows from the harmonicity of $\Delta b$ that

$$(38) \quad \lim_i \int_{\partial B_1} -\frac{\partial}{\partial \nu}(u_i(y_i)\Delta u_i) = -\int_{\partial B_1} \frac{\partial}{\partial \nu}\left(2a(n-4)|y|^{2-n} - \Delta b(y)\right) = 2a(n-4)(n-2)\omega_{n-1}.$$

On the other hand, by Lemma 2.12, we have

$$(39) \quad \lim_i \frac{n-4}{2} u_i(y_i) \int_{B_{r_i}} f_i u_i^{q_i} = \frac{n-4}{2} \lim_i f_i(y_i) l.$$

It follows from (38) and (39) that

$$a = \left(\lim_i k_i\right)^{\frac{4-n}{2}}.$$

The Proposition is established. ∎

**Remark 2.16** *As a consequence of Proposition 2.9 we have that Lemma 2.11 can be refined to*

$$\tau_i = O\left(u_i(y_i)^{-\frac{2}{n-4}}\right).$$

*To check this it is sufficient to repeat the proof of Lemma 2.11 and to use (22).*

We now state a technical Lemma, which proof is a simple consequence of Proposition 2.7, Lemma 2.11 and Proposition 2.9. We recall that $r_i = R_i \, u_i(y_i)^{-\frac{q_i-1}{4}}$.

**Lemma 2.17** *Under the hypotheses of Proposition 2.9, we have*

$$\int_{B_{r_i}} |y-y_i|^s u_i(y)^{q_i+1} = \begin{cases} u_i(y_i)^{-\frac{2s}{n-4}}\left(\int_{\mathbb{R}^n} |z|^s (1+k|z|^2)^{-n}\,dz + o(1)\right) & -n < s < n, \\ O\left(u_i(y_i)^{-\frac{2n}{n-4}} \log(u_i(y_i))\right) & s = n, \\ o\left(u_i(y_i)^{-\frac{2n}{n-4}}\right) & s > n \end{cases}$$

$$\int_{B_1\setminus B_{r_i}} |y-y_i|^s u_i(y)^{q_i+1} = \begin{cases} o\left(u_i(y_i)^{-\frac{2s}{n-4}}\right) & -n < s < n, \\ O\left(u_i(y_i)^{-\frac{2n}{n-4}} \log(u_i(y_i))\right) & s = n, \\ O\left(u_i(y_i)^{-\frac{2n}{n-4}}\right) & s > n, \end{cases}$$

*where* $k^2 = (\lim_i k_i)^2 = \frac{1}{2n(n+2)(n-2)} \lim_i f_i(y_i)$.

Now we show that under some regularity assumptions on $f_i$, $\nabla f_i(y_i)$ is small for $i$ large.

**Lemma 2.18** *Let $\overline{r} \in (0, \varrho)$, assume that $\{f_i\}_i$ is bounded in $C^2(B_{\overline{r}}(\overline{y}))$, and $u_i$ satisfy equation (3). Let $y_i \to \overline{y} \in \Omega$ be an isolated simple blow up point of $u_i$. Then*

$$|\nabla f_i(y_i)| = O\left(u_i(y_i)^{-\frac{2}{n-4}}\right).$$

PROOF  Consider a cut off function $\eta \in C_c^\infty(B_1)$, satisfying

$$\eta(x) = \begin{cases} 1 & |x| \leq \frac{1}{4}, \\ 0 & |x| \geq \frac{1}{2}. \end{cases}$$

Multiplying equation (3) by $\eta \frac{\partial u_i}{\partial x_1}$ and integrating by parts, it follows from Proposition 2.9 that

$$\frac{1}{q_i+1} \int_{B_1} \frac{\partial f_i}{\partial x_1} u_i^{q_i+1} \eta = -\int_{B_1 \setminus B_{\frac{1}{4}}} \Delta u_i \frac{\partial u_i}{\partial x_1} \Delta \eta + \frac{1}{2} \int_{B_1 \setminus B_{\frac{1}{4}}} \frac{\partial \eta}{\partial x_1} (\Delta u_i)^2$$

$$- 2 \int_{B_1 \setminus B_{\frac{1}{4}}} \Delta u_i \langle \nabla \eta, \nabla \left( \frac{\partial u_i}{\partial x_1} \right) \rangle - \frac{1}{q_i+1} \int_{B_1 \setminus B_{\frac{1}{4}}} \frac{\partial \eta}{\partial x_1} f_i u_i^{q_i+1}$$

$$= O\left( u_i(y_i)^{-2} \right).$$

Therefore, taking into account the boundedness of $\{f_i\}_i$ in $C^2(B_1)$ and Lemma 2.17 we have

$$\frac{1}{q_i+1} \int_{B_1} \eta u_i^{q_i+1} \frac{\partial f_i}{\partial x_1}(y_i) = \frac{1}{q_i+1} \int_{B_1} \left( \frac{\partial f_i}{\partial x_1}(y_i) - \frac{\partial f_i}{\partial x_1}(y) \right) \eta u_i^{q_i+1} + \frac{1}{q_i+1} \int_{B_1} \frac{\partial f_i}{\partial x_1}(y) u_i^{q_i+1} \eta$$

$$= O\left( \int_{B_1} |y - y_i| u_i(y)^{q_i+1} dy \right) + O\left( u_i(y_i)^{-2} \right) = O\left( u_i(y_i)^{-\frac{2}{n-4}} \right).$$

Clearly, we can estimate $\left( \frac{\partial f_i}{\partial x_k} \right)(y_i)$, $2 \leq k \leq n$, in a similar way, so Lemma 2.18 follows immediately. ∎

We notice that according to the above Lemma, if $f_i \to f$ in $C^2(\overline{B_{\overline{r}}(\overline{y})})$, then $\overline{y}$ is a critical point of $f$.

**Proposition 2.19** *Assume that $\{u_i\}_i$ satisfies equation (3) with $n = 5, 6$, and*

*i) if $n = 5$, $\{f_i\}_i$ is bounded in $C^1_{loc}(B_2)$ ;*

*ii) if $n = 6$, $\{f_i\}_i$ is bounded in $C^2_{loc}(B_2)$ .*

*Let $\overline{y}$ be an isolated blow up point for $\{u_i\}_i$. Then $\overline{y}$ is an isolated simple blow up point.*

PROOF  It follows from Proposition 2.7 that $r^{\frac{4}{q_i-1}} \overline{u}_i(r)$ has precisely one critical point in the interval $(0, r_i)$, where $r_i = R_i u_i(y_i)^{-\frac{q_i-1}{4}}$, as before. Arguing by contradiction, suppose that $\overline{y}$ is not an isolated simple blow up, and let $\mu_i$ be the second critical point of $r^{\frac{4}{q_i-1}} \overline{u}_i(r)$. We know that $\mu_i \geq r_i$ and, by the contradiction argument, $\mu_i \to 0$. Without loss of generality, we assume that $y_i = 0$. Set

$$\xi_i(y) = \mu_i^{\frac{4}{q_i-1}} u_i(\mu_i y) \quad |y| \leq \frac{1}{\mu_i}.$$

It follows from (3), the definition of isolated simple blow up and from the properties of $\mu_i$ that $\xi_i$ satisfies

(40) $\quad \begin{cases} \Delta^2 \xi_i(y) = f_i(\mu_i y) \xi_i(y)^{q_i} & |y| < \frac{1}{\mu_i}, \\ |y|^{\frac{4}{q_i-1}} \xi_i(y) \leq A_3 & |y| < \frac{1}{\mu_i}, \\ \lim_i \xi_i(0) = +\infty. \end{cases}$

Moreover, by our choice of $\mu_i$ there holds

$$r^{\frac{4}{q_i-1}} \overline{\xi}_i(r) \text{ has precisely one critical point in } 0 < r < 1,$$

$$\frac{d}{dr} \left( r^{\frac{4}{q_i-1}} \overline{\xi}_i(r) \right) \big|_{r=1} = 0,$$

where $\overline{\xi}_i(r) = \frac{1}{|\partial B_r|} \int_{\partial B_r} \xi$.

It follows that 0 is an isolated simple blow up for $\{\xi_i\}_i$. Therefore, applying Proposition 2.9, there exist some positive constant $a > 0$, and some regular biharmonic function $b(y)$ in $\mathbb{R}^n$ such that

(41) $\quad \xi_i(0) \xi_i(y) \to h(y) = a |x|^{4-n} + b(y) \quad \text{in } C^4_{loc}(\mathbb{R}^n \setminus \{0\}).$

We notice that $h(y)$ is positive, and $\Delta h(y)$ is non-negative.

We now claim that $b \equiv c$, for some $c \in \mathbb{R}, c \geq 0$. Indeed, we have that $\Delta b(y)$ is bounded from below by a fixed constant and is harmonic, so by the Liouville Theorem $\Delta b \equiv d$, for some $d \geq 0$. We also remark that $b(y)$ is bounded from below. Hence, if we consider the function $l(y) = b(y) + \frac{d}{2n}|y|^2$, we have that $l(y)$ is bounded from below in $\mathbb{R}^n$ and is harmonic so, again by the Liouville Theorem, it follows that $l(y)$ is constant. Since $b$ is bounded from below, this implies that $d = 0$. Hence $b(y)$ is an harmonic function in $\mathbb{R}^n$ bounded from below, and so it must be a constant.

The value of $b$ can be computed as follows. Since, by our choice of $\mu_i$, $1$ is a critical point of $r^{\frac{4}{q_i-1}}\overline{\xi}_i(r)$, we have that
$$0 = \lim_i \frac{d}{dr}\left(r^{\frac{4}{q_i-1}}\overline{\xi}_i(r)\right)|_{r=1} = \frac{d}{dr}\left(r^{\frac{n-4}{2}}\overline{h}(r)\right)|_{r=1};$$
it follows immediately that
$$b = a > 0.$$

Applying Propositions 2.9 and 4.1 to equation (40) we have, for any $\sigma \in (0, 1)$
$$\int_{\partial B_\sigma} B(\sigma, x, \xi_i, \nabla \xi_i, \nabla^2 \xi_i, \nabla^3 \xi_i) \geq \frac{n-4}{2(q_i+1)} \sum_j \int_{B_\sigma} x_j \frac{\partial f_i(\mu_i \cdot)}{\partial x_j} \xi_i^{q_i+1} - \frac{\sigma(n-4)}{2(q_i+1)} \int_{\partial B_\sigma} f_i(\mu_i \cdot) \xi_i^{q_i+1}$$
$$\geq \frac{n-4}{2(q_i+1)} \sum_j \int_{B_\sigma} x_j \frac{\partial f_i(\mu_i \cdot)}{\partial x_j} \xi_i^{q_i+1} - O\left(\xi_i(0)^{-q_i-1}\right).$$

Multiply the above by $\xi_i(0)^2$ and send $i$ to infinity, we have
$$\int_{\partial B_\sigma} B(\sigma, x, h, \nabla h, \nabla^2 h, \nabla^3 h) = \lim_i \xi_i(0)^2 \int_{\partial B_\sigma} B(\sigma, x, \xi_i, \nabla \xi_i, \nabla^2 \xi_i, \nabla^3 \xi_i)$$
$$\geq \lim_i \xi_i(0)^2 \frac{n-4}{2(q_i+1)} \sum_j \int_{B_\sigma} x_j \frac{\partial f_i(\mu_i \cdot)}{\partial x_j} \xi_i^{q_i+1}.$$

Now we want to estimate the last expression.

For $n = 5$, we recall that we are assuming that $\{f_i\}_i$ is uniformly bounded in $C^1_{loc}(B_2)$ and so, taking into account Lemma 2.17, there holds

(42) $$\left|\sum_j \int_{B_\sigma} x_j \frac{\partial f_i(\mu_i \cdot)}{\partial x_j} \xi_i^{q_i+1}\right| \leq C \mu_i \int_{B_\sigma} |x| \, \xi_i^{q_i+1} = o\left(\xi_i(0)^{-2}\right).$$

It follows that

(43) $$\int_{\partial B_\sigma} B(\sigma, x, h, \nabla h, \nabla^2 h, \nabla^3 h) \geq \lim_i \xi_i(0)^2 \frac{n-4}{2(q_i+1)} \sum_j \int_{B_\sigma} x_j \frac{\partial f_i(\mu_i \cdot)}{\partial x_j} \xi_i^{q_i+1} = 0.$$

For $n = 6$, we recall that we are assuming that $\{f_i\}_i$ is uniformly bounded in $C^2_{loc}(B_2)$; we proceed as follows. We have, using the Taylor expansion of $f_i$ at $0$
$$\left|\sum_j \int_{B_\sigma} x_j \frac{\partial f_i(\mu_i \cdot)}{\partial x_j} \xi_i^{q_i+1}\right| \leq \mu_i \left|\int_{B_\sigma} y \cdot \nabla f_i(0) \xi_i^{q_i+1}\right| + \mu_i^2 \max_{B_{\sigma\mu_i}} |\nabla^2 f_i| \int_{B_\sigma} |y|^2 \, \xi_i^{q_i+1}.$$

Applying Lemma 2.18 we have that $|\nabla f_i(0)| = O\left(\xi_i(0)^{-\frac{2}{n-4}}\right)$, therefore using Lemma 2.17 and the fact that $\mu_i \to 0$, we have

(44) $$\int_{\partial B_\sigma} B(\sigma, x, h, \nabla h, \nabla^2 h, \nabla^3 h) \geq \lim_i \xi_i(0)^2 \frac{n-4}{2(q_i+1)} \sum_j \int_{B_\sigma} x_j \frac{\partial f_i(\mu_i \cdot)}{\partial x_j} \xi_i^{q_i+1} = 0.$$

By Corollary 4.2, we know that for $\sigma > 0$ sufficiently small

$$\int_{\partial B_\sigma} B(\sigma, x, h, \nabla h, \nabla^2 h, \nabla^3 h) < 0,$$

and this contradicts both (43) and (44). This concludes the proof of the Proposition. ∎

## 3 Applications

Once the local blow up analysis is performed, we can adapt to this framework the existence arguments developed in [16] and in [30] for the scalar curvature equation. As remarked in the Introduction, the main difference with respect to the scalar curvature problem is to prove the positivity of the solutions involved in such a scheme. For $n = 5$, this is related to the estimates in [25], while for $n = 6$ this is proved in Proposition 3.6. The main ingredient of these arguments are the a priori estimates given in the next subsection.

### 3.1 A priori estimates on $\mathbb{S}^5, \mathbb{S}^6$

Consider for $n \geq 5$ the following equation

(45)
$$\begin{cases} P_h^n v - \frac{n-4}{2} f(x) v^q = 0 & \text{on } \mathbb{S}^n, \\ v > 0, & \text{on } \mathbb{S}^n, \\ q = \frac{n+4}{n-4} - \tau, \quad 0 \leq \tau < \frac{8}{n-4}. \end{cases}$$

**Proposition 3.1** *Suppose $f \in C^1(\mathbb{S}^n)$ satisfies for some positive constant $A_1$*

$$\frac{1}{A_1} \leq f(p) \leq A_1 \quad \text{for all } p \in \mathbb{S}^n.$$

*Then, for any $0 < \varepsilon < 1$, $R > 1$, there exist some positive constants $C_0^* > 1$, $C_1^* > 1$ depending on $n, \varepsilon, R, A_1, \|f\|_{C^1}$ such that, if $v$ is a solution of (45) with*

$$\max_{\mathbb{S}^n} v > C_0^*,$$

*then there exists $1 \leq k = k(v) < +\infty$ and a set*

$$S(v) = \{p_1, \ldots, p_k\} \subseteq \mathbb{S}^n, \qquad p_i = p_i(v)$$

*such that*

(1) $0 \leq \tau < \varepsilon$,

(2) $p_1, \ldots, p_k$ *are local maxima of $v$ and for each $1 \leq j \leq k$, using $y$ as some geodesic normal coordinates centered at $p_j$, we have*

$$\left\| v(0)^{-1} v\left(v(0)^{-\frac{q-1}{2}} y\right) - \delta_j(y) \right\|_{C^4(B_{2R}(0))} < \varepsilon$$

*and*

$$\left\{ B_{R\, v(p_j)^{-\frac{q-1}{4}}}(p_j) \right\}_{1 \leq j \leq k} \quad \text{are disjoint balls.}$$

*Here*
$$\delta_j(y) = \left(1 + k_j |y|^2\right)^{\frac{4-n}{2}}$$

*is the unique solution of*

$$\begin{cases} \Delta^2 \delta_j = \frac{n-4}{2} f(p_i) \delta_j^{\frac{n+4}{n-4}} & \text{in } \mathbb{R}^n, \\ \delta_j > 0, \quad \delta_j(0) = 1, \quad \nabla \delta_j(0) = 0 & \text{in } \mathbb{R}^n, \end{cases}$$

*and*

$$k_j^2 = \frac{1}{2n(n+2)(n-2)} f(p_j),$$

(3) $v(p) \leq C_1^* \{dist(p, \mathcal{S}(v))\}^{-\frac{4}{q_i-1}}$ *for all* $p \in \mathbb{S}^n$, *and* $dist(p_i, p_j)^{\frac{4}{q_i-1}} v(p_j) \geq C_0^*$.

PROOF  This can be proved by quite standard blow up arguments, using the results of [31], see [38]. ∎

**Proposition 3.2** *Let $n = 5, 6$ and assume that $f \in C^1(\mathbb{S}^n)$ satisfies, for some positive constant $A_1$*

$$\frac{1}{A_1} \leq f(p) \leq A_1 \quad \forall p \in \mathbb{S}^n.$$

*If $n = 6$, we also assume that $f$ is of class $C^2$ on $\mathbb{S}^n$. Then for every $\varepsilon > 0$ and $R > 1$, there exists some positive constant $\delta^* > 0$, depending on $n, \varepsilon, R$ and $\|f\|_{C^1}$ (and also on $\|f\|_{C^2}$ if $n = 6$), such that for any solution $v$ of (45) with $\max_{\mathbb{S}^n} v > C_0^*$ we have*

$$|p_j - p_l| \geq \delta^* \quad \text{for all } 1 \leq j \neq l \leq k,$$

*where $C_0^*$, $p_j = p_j(v), p_l = p_l(v)$ and $k = k(v)$ are as in Proposition 3.1.*

PROOF  Suppose the contrary, that is for some constants $\varepsilon, R, A_1$ there exist $\{q_i\}, \{f_i\}$ satisfying the assumptions of Proposition 3.2 and a sequence of corresponding solutions $v_i$ such that

$$\lim_i \min_{j \neq l} |p_j - p_l| = 0.$$

Without loss of generality, we assume that

(46) $$|p_1(v_i) - p_2(v_i)| = \min_{j \neq l} |p_j(v_i) - p_l(v_i)| \to 0.$$

Since $B_{R v_i(p_1)^{-\frac{q_i-1}{4}}}(p_1)$ and $B_{R v_i(p_2)^{-\frac{q_i-1}{4}}}(p_2)$ are disjoint we have by (46) $v_i(p_1) \to +\infty$ and $v_i(p_2) \to +\infty$.

Performing a stereographic projection with $p_1$ as the south pole and using formula (7), equation (45) is transformed into

(47) $$\begin{cases} \Delta^2 u_i = \frac{n-4}{2} f(x) H(x)^{\tau_i} u_i^{q_i} & \text{on } \mathbb{R}^n, \\ v > 0 & \text{on } \mathbb{R}^n, \\ q_i = \frac{n+4}{n-4} - \tau_i, \quad 0 \leq \tau \leq \frac{2}{n-4}. \end{cases}$$

Let us still use $p_2 \in \mathbb{R}^n$ to denote the stereographic coordinates of $p_2 \in \mathbb{S}^n$, and set $\sigma_i = |p_2| \to 0$. For simplicity we assume that $p_2$ is a local maximum of $u_i$. It is easy to see that

(48) $$\sigma_i > \frac{1}{C(n)} \max \left\{ R u_i(0)^{-\frac{q_i-1}{4}}, R u_i(p_2)^{-\frac{q_i-1}{4}} \right\}.$$

Set now
$$w_i(y) = \sigma_i^{\frac{4}{q_i-1}} u_i(\sigma_i y) \quad |y| < \frac{1}{\sigma_i}.$$

It follows that $w_i$ satisfies

(49)
$$\begin{cases} \Delta^2 w_i = \frac{n-4}{2} f(\sigma_i y) H(\sigma_i y)^{\tau_i} w_i^{q_i} & |y| \leq \frac{1}{\sigma_i}, \\ w_i > 0, \quad \Delta w_i > 0, & |y| \leq \frac{1}{\sigma_i}. \end{cases}$$

Notice that by Proposition 3.1 there holds

(50) $$v_i(y) \leq C_1 |y|^{-\frac{4}{q_i-1}} \quad \text{for all } |y| \leq \frac{1}{2}\sigma_i,$$

(51) $$v_i(y) \leq C_1 |y - p_2|^{-\frac{4}{q_i-1}} \quad \text{for all } |y - p_2| \leq \frac{1}{2}\sigma_i.$$

It is not difficult to see, using (48) and the last estimates, that
$$w_i(0) \geq C_2, \quad \lim_i w_i(|p_2|^{-1} p_2) \geq C_2,$$
$$|y|^{\frac{4}{q_i-1}} w_i(y) \leq C_1 \quad |y| \leq \frac{1}{2},$$
$$|y - |p_2|^{-1} p_2|^{\frac{4}{q_i-1}} w_i(y) \leq C_1 \quad |y - |p_2|^{-1} p_2| \leq \frac{1}{2}.$$

We first show that

(52) $$w_i(0) \to \infty, \quad w_i(p_2 |p_2|^{-1}) \to \infty.$$

If one of these sequences tends to infinity along a subsequence, say $w_i(0) \to \infty$, then 0 is an isolated simple blow up. Therefore $w_i(|p_2|^{-1}p_2)$ must tend to infinity along the same subsequence, since otherwise, using (50), it is easy to prove that $(w_i)$ is uniformly bounded near $|p_2|^{-1}p_2$ along a further subsequence. In turn, using Proposition 2.9 and the Harnack Inequality we obtain that $w_i$ tends to 0 uniformly near $|p_2|^{-1}p_2$, which is impossible. On the other hand if both $w_i(0)$ and $w_i(|p_2|^{-1}p_2)$ stay bounded, $w_i$ will be locally bounded and will converge in $C^2_{loc}$ to some function having at least two critical points, which violates the uniqueness result of C. S. Lin [31]. We thus established (52).

It then follows from Proposition 2.19 that 0 and $\overline{q} = \lim_i |p_2|^{-1} p_2$ are both isolated blow up points for $w_i$. Hence by Proposition 2.19 they are indeed isolated simple blow up points of $w_i$.
We deduce from property (3) in Proposition 3.1, (46), Lemma 2.5 and Proposition 2.9 that there exist an at most countable set $\mathcal{S}_1 \subseteq \mathbb{R}^n$ such that
$$\min\{|x - y| : x, y \in \mathcal{S}_1\} \geq 1,$$
$$\lim_i w_i(0) w_i(y) = h^*(y) \quad \text{in } C^4_{loc}(\mathbb{R}^n \setminus \mathcal{S}_1),$$
$$h^*(y) \geq 0, \quad \Delta h^*(y) \geq 0, \quad \Delta^2 h^*(y) = 0 \quad y \in \mathbb{R}^n \setminus \mathcal{S}_1,$$

and $h^*(y)$ must be singular at 0 and at $\overline{p} = \lim_i |p_2|^{-1} p_2$ ($0, \overline{p} \in \mathcal{S}_1$). Hence for some positive constants $a_1$ and $a_2$ it is

(53) $$h^*(y) = a_1 |y|^{4-n} + a_2 |y - \overline{q}|^{4-n} + b^*(y) \quad y \in \mathbb{R}^n \setminus \{\mathcal{S}_1 \setminus \{0, \overline{p}\}\},$$

where $b^*(y) : \mathbb{R}^n \setminus \{\mathcal{S}_1 \setminus \{0, \overline{p}\}\}$ satisfies
$$\Delta^2 b^*(y) = 0, \quad \liminf_{|y| \to +\infty} b^*(y) \geq 0, \quad \liminf_{|y| \to +\infty} \Delta b^*(y) \geq 0.$$

The maximum principle, applied first to $\Delta b^*(y)$ and then to $b^*(y)$, yields $b^*(y) \geq 0$ in $\mathbb{R}^n \setminus \{\mathcal{S}_1 \setminus \{0, \overline{p}\}\}$. Hence the conclusion follows from (53), reasoning as in the proof of Proposition 2.19. ∎

PROOF OF THEOREM 1.5 Proposition 3.2 and statement (3) in Proposition 3.1 imply that the blow up points are in finite number and are isolated. Hence by Proposition 2.19 they are also isolated simple. Then the conclusion follows from Proposition 2.7 and Lemma 2.17. ∎

## 3.2 Existence and compactness results on $\mathbb{S}^5$

This subsection is devoted to the proof of our existence and compactness results on $S^5$. The first ingredient is the fact that there is at most one blow up point; this is the content of Theorem 1.6.

PROOF OF THEOREM 1.6 Assume the contrary, namely that $\{v_i\}$ has $p^1, p^2 \in \mathbb{S}^n$, $p^1 \neq p^2$ as blow up points. Let $p_i^1 \to p^1$, $p_i^2 \to p^2$ be local maxima of $v_i$ as before. Without loss of generality, we may assume that $p^2 \neq -p^1$. We make a stereographic projection with $p_i^1$ being the south pole. In the stereographic coordinates, it is $p_i^1 = 0$ and we still use the notation $p_i^2, p^2$ for the projection of those points. According to (7), equation (6) becomes

$$\Delta^2 u_i(y) = \frac{n-4}{2} f_i(y) H_i(y)^{\tau_i} u_i(y)^{q_i} \quad \text{in } \mathbb{R}^n.$$

As before we can assume that $p_i^2$ is a local maximum of $u_i$. We recall that, by Proposition 3.2, the number of blow up points is bounded by some constant independent of $i$. Therefore, reasoning as above, there exist some finite set $\mathcal{S}_2 \subseteq \mathbb{R}^n$, $0, p^2 \in \mathcal{S}_2$, some constants $a_1, A > 0$ and some function $h(y) \in C^4(\mathbb{R}^n \setminus \mathcal{S}_2)$ such that

$$\lim_i u_i(0) u_i(y) = h(y) \quad \text{in } C^4_{loc}(\mathbb{R}^n \setminus \mathcal{S}_2),$$

$$h(y) = a_1 |y|^{4-n} + A + O(|y|) \quad \text{for } y \text{ close to } 0.$$

It follows from the proof of Proposition 2.19 that for any $0 < \sigma < 1$, we have

$$\lim_i u_i(0)^2 \int_{B_\sigma} y \cdot \nabla(f_i H_i^{\tau_i}) u_i^{q_i+1} = 0.$$

From Corollary 4.2 we obtain a contradiction as before. Theorem 1.6 is therefore established. ∎

**Theorem 3.3** *Let $n = 5, 6$, and assume that $f \in C^2(\mathbb{S}^n)$ is a positive Morse function which satisfies (ND). Let $\{f_i\}_i$ be a sequence of Morse functions converging to $f$ in $C^2(\mathbb{S}^n)$, and let $v_i$ be a solution of (P) with $f = f_i$. Then, after passing to a subsequence, either $v_i$ stays bounded in $L^\infty(\mathbb{S}^n)$, or has at least two blow up points.*

PROOF Suppose by contradiction that $v_i$ has precisely one blow up point $q_0$. Making a stereographic projection with $q_0$ being the south pole, the equation is then transformed into

$$(54) \qquad \Delta^2 u_i = \frac{n-4}{2} f_i(y) u_i(y)^{\frac{n+4}{n-4}}, \quad u_i > 0, \quad \Delta u_i > 0 \quad \text{in } \mathbb{R}^n.$$

Here we have identified $f_i$ with its composition with the stereographic projection. We know by Theorem 1.6 that $y_i \to 0$ is an isolated simple blow up point for $u_i$. It follows from Lemma 2.18 that

$$|\nabla f_i(y_i)| = O\left(u_i(y_i)^{-\frac{2}{n-4}}\right).$$

We can suppose without loss of generality that $q_0$ is a critical point of $f_i$: hence, from the uniform invertibility of $\nabla^2 f_i$ in $q_0$ we deduce that

$$(55) \qquad |y_i| = O\left(u_i(y_i)^{-\frac{2}{n-4}}\right).$$

Let $\xi = \lim_i u_i(y_i)^{\frac{2}{n-4}} y_i$, and let $Q(x) = (D^2 f(0) x, x)$. By means of (55), following the computations in [30], pages 370-373, we deduce that $\xi$ and $Q$ satisfy

$$(56) \qquad \int_{\mathbb{R}^n} \nabla Q(x + y_i) (1 + k|z|^2)^{-n} = 0.$$

(57) $$\int_{\mathbb{R}^n} (z+\xi)\nabla Q(x+y_i)\left(1+k\left|z\right|^2\right)^{-n} = 0.$$

It is easy to see that (56) and (57) contradict the hypothesis $(ND)$. ∎

PROOF OF THEOREM 1.7 It is an immediate consequence of Theorems 1.6 and 3.3. ∎

PROOF OF THEOREM 1.8 For $\mu \in [0,1]$, consider the function $f_\mu$ defined by

$$f_\mu = \mu f + (1-\mu)\frac{105}{8}.$$

Let $\alpha \in (0,1)$ be fixed. By Theorem 1.7, for every $\mu_0 \in ]0,1[$ there exists a positive constant $C(\mu_0)$ such that every positive solution $v$ of $(P)$ with $f = f_\mu$ and $\mu \geq \mu_0$ satisfies

(58) $$\frac{1}{C(\mu_0)} \leq v \leq C(\mu_0), \quad \|v\|_{C^{4,\alpha}(\mathbb{S}^5)} \leq C(\mu_0).$$

Using the $L^\infty$ estimates in [25], we can follow the arguments in [16] (Section 7) to check that for $\mu_0$ sufficiently small there holds

$$\deg_{C^{4,\alpha}}\left(u - (P_h^5)^{-1}(\frac{1}{2}f_{\mu_0}\left|u\right|^8 u), \left\{\frac{1}{C(\mu_0)} \leq v \leq C(\mu_0)\right\}, 0\right) = (-1)^n \deg\left(\int_{\mathbb{S}^n} f \circ \varphi_{P,t}(x)\, x, B, 0\right).$$

Under assumption (8), it follows that $\deg(\int_{\mathbb{S}^n} f \circ \varphi_{P,t}(x)\, x, B, 0) \neq 0$, see [16].

By Theorem 1.7, $u - (P_h^5)^{-1}(\frac{1}{2} f_\mu \left|u\right|^8 u)$ is different from zero on the boundary of $\left\{\frac{1}{C(\mu)} \leq v \leq C(\mu)\right\}$ hence, from the homotopy property of the degree, we have also

$$\deg_{C^{4,\alpha}}\left(u - (P_h^5)^{-1}(\frac{1}{2} f \left|u\right|^8 u), \left\{\frac{1}{C(\mu)} \leq v \leq C(\mu)\right\}, 0\right) \neq 0.$$

This concludes the proof. ∎

### 3.3 Existence and compactness results on $\mathbb{S}^6$

This subsection is devoted to prove the existence and compactness results on $\mathbb{S}^6$. Similarly to the four dimensional case for scalar curvature, there could be solutions of (6) blowing up at more than one point. In the following Proposition, we give necessary conditions for solutions to blow up, and we locate their blow up points.

**Proposition 3.4** *Let $f \in C^2(\mathbb{S}^6)$ be a positive function. Then there exists some number $\delta^* > 0$, depending only on $\min_{\mathbb{S}^6} f$ and $\|f\|_{C^2(\mathbb{S}^6)}$, with the following properties.*
*Let $\{q_i\}$ satisfy $q_i \leq 5$, $q_i \to 5$, $\{f_i\}_i \in C^2(\mathbb{S}^6)$ satisfy $f_i \to f$ in $C^2(\mathbb{S}^6)$, $v_i$ satisfy*

(59) $$P_h^n v_i = f\, v_i^{q_i}, \quad v_i > 0 \quad \text{on } \mathbb{S}^6,$$

*and $\limsup_i \max_{\mathbb{S}^6} v_i = +\infty$. Then after passing to a subsequence, we have*

i) *$\{v_i\}_i$ has only isolated simple blow up points $(p^1, \ldots, p^k) \in \mathcal{F}\setminus\mathcal{F}^-$ $(k \geq 1)$, with $|p^j - p^i| \geq \delta^*$ $\forall j \neq k$, and $\rho(p^1, \ldots, p^k) \geq 0$. Furthermore $p^1, \ldots, p^k \in \mathcal{F}^+$ if $k \geq 2$.*

ii) *Setting*
$$\lambda_j = f(p^j)^{-\frac{1}{4}} \lim_i v_i(p_i^1)\,(v_i(p_i^j))^{-1}, \quad \mu^j := \lim_i \tau_i\, v_i(p_i^j)^2,$$

*where $p_i^j \to p^j$ is the local maximum of $v_i$, there holds*

$$\lambda_j \in\, ]0, +\infty[, \quad \mu^j \in [0, +\infty[ \quad \forall j = 1, \ldots, k.$$

*iii)* When $k = 1$
$$\mu^1 = \frac{\Delta f(p^1)}{(f(p^1))^{\frac{3}{2}}};$$

when $k \geq 2$

(60) $$\sum_{l=1}^{k} M_{lj}\lambda_l = 4\frac{\sqrt{6}}{3}\lambda_j\mu^j \quad \forall j = 1, \ldots, k.$$

*iv)* $\mu^j \in ]0, +\infty[ \ \forall j = 1, \ldots, k$ if and only if $\rho(p^1, \ldots, p^k) > 0$.

PROOF Assertion *ii)* follows from Proposition 2.9, Lemma 2.5 and Remark 2.16. From another part, it follows from Proposition 3.2, Proposition 2.19 that $v_i$ has only isolated simple blow up points $p^1, \ldots, p^k \in \mathcal{F}$ ($k \geq 1$) with $|p^j - p^l| \geq \delta^*$ ($j \neq l$) for a fixed $\delta^* > 0$.

Let $p_i^1 \to p^1$ be the local maximum of $v_i$ for which $v_i(p_i^1) \to +\infty$. Making a stereographic projection with south pole $p_i^1$, equation (59) is transformed into

$$\Delta^2 u_i(y) = f_i(y) H^{\tau_i}(y) u_i^{q_i} \quad y \in \mathbb{R}^6.$$

By our choice of the projection, 0 is a local maximum for all $u_i$; moreover, it is clear that 0 is also an isolated simple blow up point. We can also suppose that none of the points $\{p^1, \ldots, p^k\}$ is mapped to $+\infty$ by the stereographic projection, and we still denote their images by $p^1, \ldots, p^k$. It follows from Proposition 2.9 that

(61) $$u_i(y_i^j) u_i(y) \to h^j(y) := 8\sqrt{6} f(p^j)^{-\frac{1}{2}} |y|^{-2} + b^j(y) \quad \text{in } C^4_{loc}(\mathbb{R}^6 \setminus \{p^1, \ldots, p^k\}),$$

where $b^j$ is some biharmonic function in $\mathbb{R}^6 \setminus \{p^2, \ldots, p^k\}$.
Coming back to $v_i$ we have

$$\lim_i v_i(p_i^1) v_i(p) = 4\sqrt{6} f(p^j)^{-\frac{1}{2}} J_{p^1}(p) + \tilde{b}^1(p) \quad \text{in } C^4_{loc}(\mathbb{S}^6 \setminus \{p^2, \ldots, p^k\}),$$

where $\tilde{b}^1$ is some regular function on $\mathbb{S}^6 \setminus \{p^2, \ldots, p^k\}$ satisfying $P_h^6 \tilde{b}^j = 0$.
If $k = 1$, then $\tilde{b}^1 = 0$ while for $k \geq 2$, taking into account the contribution of all the poles, we deduce that for all $j = 1, \ldots, k$ it is

$$\lim_i v_i(p_i^j) v_i(p) = 4\sqrt{6} \left\{ \frac{J_{p^j}(p)}{\sqrt{f(p^j)}} + \sum_{l \neq j} \lim_i \frac{v_i(p_i^j)}{v_i(p_i^l)} \frac{J_{p^l}(p)}{\sqrt{f(p^l)}} \right\},$$

where the convergence is in $C^4_{loc}(\mathbb{S}^6 \setminus \{p^1, \ldots, p^k\})$. In fact, subtracting all the poles from the limit function, we obtain a regular function $r: \mathbb{S}^6 \to \mathbb{R}$ for which $P_h^6 r = 0$; by the coercivity of $P_h^6$ on $H_2^2(\mathbb{S}^6)$ it must be $r = 0$.
Using the last formula, we can compute the exact expression of $h^j(y)$, which is

(62) $$h^j(y) = 8\sqrt{6} f(p^j)^{-\frac{1}{2}} |y|^{-2} + 16\sqrt{6} \sum_{l \neq j} \lim_i \frac{v_i(p_i^j)}{v_i(p_i^l)} \frac{J_{p_l}(p^j)}{\sqrt{f(p^l)}} + O(|y|).$$

Hence, using (62) and Corollary 4.2, we deduce that

$$\lim_{\sigma \to 0} \int_{\partial B_\sigma} B(\sigma, x, h^j, \nabla h^j, \nabla^2 h^j, \nabla^3 h^j) = -3 \, (2^9) \, \omega_5 \sum_{l \neq j} \lim_i \frac{v_i(p_i^j)}{v_i(p_i^l)} \frac{J_{p_l}(p^j)}{\sqrt{f(p^l)}\sqrt{f(p^j)}}.$$

From another part, it follows from Propositions 4.1, 2.7 and Lemma 2.11 that, for any $0 < \sigma < 1$

$$\int_{\partial B_\sigma} B(\sigma, x, h^j, \nabla h^j, \nabla^2 h^j, \nabla^3 h^j) = -\frac{2^9}{5} \omega_5 \frac{\Delta_h f(p^j)}{f(p^j)^2} + \frac{\sqrt{6}\, 2^{11}}{15} \omega_5 \frac{\mu^j}{\sqrt{f(p^j)}}.$$

From the last two formulas, using the expression of $\mu_l$ and $\lambda_l$, we obtain

$$-15 \left( \sum_{l \neq j} \frac{J_{p^l}(p^j)}{f(p^j)^{\frac{1}{4}} f(p^l)^{\frac{1}{4}}} \lambda_l \right) + \frac{\Delta_h f(p^j)}{f(p^j)^{\frac{3}{2}}} \lambda_j = 4\frac{\sqrt{6}}{3} \mu^j \lambda_j.$$

We have thus established (60); in particular when $k = 1$ we obtain $\mu^1 = \frac{\Delta_h f(p^1)}{f(p^1)^{\frac{3}{2}}}$, so we have deduced $iii$).

It follows that $p^j \in \mathcal{F} \setminus \mathcal{F}^-$, $\forall j = 1, \ldots, k$, and when $k \geq 2$, $p^j \in \mathcal{F}^+$. Furthermore, since $M_{ii} \geq 0$, and $M_{ij} < 0$ for $i \neq j$, it follows from linear algebra and the variational characterization of the least eigenvalue that there exists some $x = (x_1, \ldots, x_k) \neq 0$, $x_l \geq 0$ $\forall l$, such that $\sum_{j=1}^k M_{lj} x_j = \rho\, x_l$. Multiplying (60) by $x_j$ and summing over $j$, we have

$$\rho \sum_j \lambda_j x_j = \sum_{l,j} M_{lj} \lambda_j x_l = 2\frac{\sqrt{6}}{3} \sum_j \lambda_j x_j \mu^j \geq 0.$$

It follows that $\rho \geq 0$, so we have verified part $i$). Part $iv$) follows from $i$)-$iii$). ■

Now we perform the following construction, needed in the proof of Theorem 1.9. For $a \in \mathbb{S}^6$ and $\lambda > 0$, let $\varphi_{a,\lambda} : \mathbb{S}^6 \to \mathbb{S}^6$ be the conformal transformation defined in the introduction, and let

$$\delta_{a,\lambda}(x) = |\det d\varphi_{a,\lambda^{-1}}|^{\frac{1}{6}}.$$

For all the choices of $a$ and $\lambda$, the function $\delta_{a,\lambda}$ satisfies $P_h^6 \delta_{a,\lambda} = 24\, \delta_{a,\lambda}^5$. We consider the following scalar product and norm on $H_2^2(\mathbb{S}^6)$ which is equivalent to the usual one, see [24],

$$\langle u, v \rangle = \int_{\mathbb{S}^6} (P_h^n u)\, v, \quad \|u\| = \langle u, u \rangle^{\frac{1}{2}}.$$

Set for $\tau > 0$ small

$$I_\tau(u) = \frac{1}{2} \int_{\mathbb{S}^6} (\Delta u)^2 + 10\, |\nabla u|^2 + 24\, u^2 - \frac{1}{6-\tau} \int_{\mathbb{S}^6} f\, |u|^{6-\tau} \quad u \in H_2^2(\mathbb{S}^6).$$

Let $p^1, \ldots, p^k \in \mathcal{F}^+$ be critical points of $f$ with $\rho(p^1, \ldots, p^k) > 0$. For $\varepsilon_0$ small, let $V_{\varepsilon_0} = V_{\varepsilon_0}(p^1, \ldots, p^k) \in \mathbb{R}_+^k \times \mathbb{R}_+^k \times (\mathbb{S}^6)^k$ be defined by

$$V_{\varepsilon_0} = \left\{ (\alpha, \lambda, a) \in \mathbb{R}_+^k \times \mathbb{R}_+^k \times (\mathbb{S}^6)^k : \left| \alpha_i - \left(\frac{24}{f(a_i)}\right)^{\frac{1}{4}} \right| < \varepsilon_0, |a_i - p^i| < \varepsilon_0, \lambda_i > \frac{1}{\varepsilon_0}, i = 1, \ldots, k \right\}.$$

It follows arguing as in [8], [10], that there exists $\varepsilon_0 > 0$ small, depending only on $\min_{\mathbb{S}^6} f$, and $\|f\|_{C^2(\mathbb{S}^6)}$, with the following property. For any $u \in H_2^2(\mathbb{S}^6)$ satisfying for some $(\tilde{\alpha}, \tilde{\lambda}, \tilde{a}) \in V_{\frac{\varepsilon_0}{2}}$ the inequality $\left\| u - \sum_{i=1}^k \tilde{\alpha}_i\, \delta_{\tilde{a}_i, \tilde{\lambda}_i} \right\| < \frac{\varepsilon_0}{2}$, we have a unique representation

$$u = \sum_{i=1}^k \alpha_i\, \delta_{a_i, \lambda_i} + v,$$

with $(\alpha, \lambda, a) \in V_{\varepsilon_0}$ and

$$(63) \qquad \langle v, \delta_{a_i, \lambda_i} \rangle = \langle v, \frac{\partial \delta_{a_i, \lambda_i}}{\partial a_i} \rangle = \langle v, \frac{\partial \delta_{a_i, \lambda_i}}{\partial \lambda_i} \rangle = 0.$$

Denote by $E_{\lambda, a}$ the set of $v \in H_2^2(\mathbb{S}^6)$ satisfying (63). It follows that in a small neighborhood (independent of $\tau$) of $\left\{ \sum_{i=1}^k \alpha_i \delta_{a_i, \lambda_i} : (\alpha, \lambda, a) \in \Omega_{\frac{\varepsilon_0}{2}} \right\}$, $(\alpha, \lambda, a, v)$ is a good parametrization of $u$. For a large constant $A$ and for a small constant $v_0$, set

$$\begin{aligned}
\Sigma_\tau &= \Sigma_\tau(p_1, \ldots, p_k) = \{(\alpha, \lambda, a, v) \in V_{\frac{\varepsilon_0}{2}} \times H_2^2(\mathbb{S}^6) : \\
&\quad |a_i - p_i| < \sqrt{\tau} |\log \tau|, A^{-1} \sqrt{\tau} < \lambda_i^{-1} < A \sqrt{\tau}, v \in E_{\lambda, a}, \|v\| < v_0 \}.
\end{aligned}$$

Without confusion, we use the same notation for

$$\Sigma_\tau = \left\{ u = \sum_{i=1}^k \alpha_i \delta_{a_i, \lambda_i} + v : (\alpha, \lambda, a, v) \in \Sigma_\tau \right\} \subseteq H_2^2(\mathbb{S}^6).$$

From Proposition 3.4 and Remark 2.4, one can easily deduce the following Proposition. We recall that we have set

$$\mathcal{O}_R = \left\{ v \in C^{4, \alpha}(\mathbb{S}^6) : \frac{1}{R} \leq v \leq R, \|v\|_{C^{4, \alpha}(\mathbb{S}^6)} \leq R \right\}.$$

**Proposition 3.5** *For $f \in \mathcal{A}$, $\alpha \in ]0, 1[$, there exist some positive constants $v_0 \ll 1$, $A \gg 1$, $R \gg 1$, depending only on $f$ such that when $\tau > 0$ is sufficiently small,*

$$u \in \mathcal{O}_R \cup \left\{ \Sigma_\tau(p^1, \ldots, p^k) : p^1, \ldots, p^k \in \mathcal{F}^+, \rho(p^1, \ldots p^k) > 0, k \geq 1 \right\},$$

*for all $u$ satisfying $u \in H^2(\mathbb{S}^6)$, $u > 0$ and $I_\tau'(u) = 0$.*

If $f \in \mathcal{A}$, we can also give sufficient conditions for the existence of positive solutions of $I_\tau' = 0$.

**Proposition 3.6** *Let $f \in \mathcal{A}$, $v_0 > 0$ be suitably small and $A > 0$ be suitably large. Then, if $p^1, \ldots, p^k \in \mathcal{F}^+$ with $\rho(p^1, \ldots, p^k) > 0$, and if $\tau > 0$ is sufficiently small, the functional $I_\tau$ has a unique critical point $u$ in $\Sigma_\tau$. In the above parametrization, we have $v \to 0$ as $\tau \to 0$.*
*This function $u$ is positive, and as a critical point of $I_\tau$ it is nondegenerate with Morse index $7k - \sum_{j=1}^k m(f, p^j)$, where $m(f, p^j)$ is the Morse index of $K$ at $p^j$.*

PROOF The proof of the existence and uniqueness of a non degenerate critical point is based only on the study of $I_\tau$ is $\Sigma_\tau$ and this can be performed as in [30], see also [11], so we omit it here. We just remark that it uses a local inversion theorem, which can be applied by the properties of the spectrum of the conformal laplacian on $\mathbb{S}^n$. Since the spectrum of $P_h^n$ possesses analogous properties, see [25] Theorem 2.2, we are indeed in the same situation from the variational point of view.
Differently from the scalar curvature case, the proof of the positivity is more involved, and we perform it in Subsection 3.4. This difficulty arises from the fact that we cannot use as a test function the negative part of $u$. ∎

When the number $\tau$ is bounded from below, we have also the following compactness result for positive solutions.

**Proposition 3.7** *Let $f \in C^2(\mathbb{S}^6)$ be a positive function, $0 < \tau_0 < \tau \leq 4 - \tau_0$. There exist some positive constants $C$ and $\delta$ depending only on $\tau_0$, $\min_{\mathbb{S}^6} f$, and $\|f\|_{C^2}$ with the following properties*

i) $\{u \in H_2^2(\mathbb{S}^6) : u \geq 0 \text{ a.e.}, I'_\tau(u) = 0\} \subseteq \mathcal{O}_C$,

ii) setting $\mathcal{O}_{C,\delta} = \{u \in H_2^2(\mathbb{S}^6) : \exists v \in \mathcal{O}_C \text{ such that } \|u - v\|_{H_2^2} < \delta\}$, it is $I'_\tau \neq 0$ on $\partial \mathcal{O}_{C,\delta}$, and

$$\deg_{H_2^2}(u - (P_h^6)^{-1}(f|u|^{4-\tau}u), \mathcal{O}_{C,\delta}, 0) = -1. \tag{64}$$

PROOF  Property i) is a consequence of the nonexistence results of [31] and of Remark 2.4. The fact that $I'_\tau \neq 0$ on $\partial \mathcal{O}_{C,\delta}$ is a consequence of the $L^\infty$ estimates in [25], see e.g. Lemma 4.9 there.
In fact, having uniform estimates from above and from below on the positive solutions of $I'_\tau = 0$, it is possible to prove (subtracting the equations) that solutions $u$ of $P_h^6 u = f|u|^{4-\tau}u$ which are close in $H_2^2$ to elements of $\mathcal{O}_C$ are also $L^\infty$ close. Hence they are positive and still contained in $\mathcal{O}_C$.
About the computation of the degree, consider the homotopy $f_t = tf + (1-t)f^*$, with $f^* = x^7 + 2$, recall that $\mathbb{S}^6 = \{x \in \mathbb{R}^7 : \|x\| = 1\}$. It follows from the Kazdan-Warner condition, see the Introduction, that there is no solution of $(P)$ with $f = f^*$. Therefore we only need to establish (64) for $f^*$ and $\tau$ very small. This follows from Propositions 3.4, 3.5 and 3.6. ∎

PROOF OF THEOREM 1.9 The norm inequality in (11) follows from Theorem 1.5. Suppose by contradiction that the second inequality is not true; then there exist solutions $v_i$ blowing up at $p^1, \ldots, p^k \in \mathbb{S}^6$, and these are isolated simple blow up points. It follows from Theorem 3.3 that $k \geq 2$: taking into account that $f \in \mathcal{A}$ and $\mu^j = 0$ for all $j$ ($\tau_i = 0$), we get a contradiction by Proposition 3.4 iv). Hence (11) is proved.
Using Proposition 3.5, (11) and the homotopy invariance of the Leray-Schauder degree, we have

$$\deg(u - (P_h^n)^{-1}(f|u|^4 u), \mathcal{O}_R, 0) = \deg(u - (P_h^n)^{-1}(f|u|^{4-\tau}u), \mathcal{O}_R, 0). \tag{65}$$

By Propositions 3.5 and 3.6, for suitable values of $\tau, A$ and $v_0$ we know that the positive solutions of $I'_\tau = 0$ are either in $\mathcal{O}_R$ or in some $\Sigma_\tau$, and viceversa for all $p^1, \ldots, p^k \in \mathcal{F}_+$ with $\rho(p^1, \ldots, p^k) > 0$, there is a nondegenerate critical point of $I_\tau$ in $\Sigma_\tau$ which is a positive function. This gives a complete characterization of the positive solutions of (45) when $\tau$ is positive and small.
Let $C$ and $\delta$ be given by Proposition 3.7. It is clear that if $C$ is sufficiently large and $\delta_1$ is sufficiently small, then $\mathcal{O}_{R,\delta_1} \subseteq \mathcal{O}_{C,\delta}$. By Proposition 3.6, (64) and by the excision property of the degree, we have

$$\deg_{H_2^2}(u - (P_h^n)^{-1}(f|u|^{4-\tau}u), \mathcal{O}_{R,\delta_1}, 0) = \text{Index}(f). \tag{66}$$

As in the proof of Proposition 3.7, one can check that there are no critical points of $I_\tau$ in $\overline{\mathcal{O}_{R,\delta_1}} \setminus \mathcal{O}_R$, hence Theorem B.2 of [30] Part I applies and yields

$$\deg_{H_2^2}(u - (P_h^n)^{-1}(f|u|^{4-\tau}u), \mathcal{O}_{R,\delta_1}, 0) = \deg(u - (P_h^n)^{-1}(f|u|^{4-\tau}u), \mathcal{O}_R, 0). \tag{67}$$

Then the conclusion follows from (65), (66) and (67). The proof of Theorem 1.9 is thereby completed. ∎

### 3.4  Positivity of the solutions

In this subsection we prove the positivity statement in Proposition 3.6. We define the operator $L_h$ to be $L_h u = \Delta u + \frac{c_n}{2} u$, and we consider the problem

$$L_h^2 u = g \quad \text{on } \mathbb{S}^n, \tag{68}$$

where $g \in L^p(\mathbb{S}^n)$, for some $p > 1$. From standard elliptic theory there exists an unique weak solution $u \in H_4^p(\mathbb{S}^n)$, and moreover

$$\|u\|_{H_4^p} \leq \overline{C}(n,p) \|g\|_p, \tag{69}$$

for some constant $\overline{C}(n,p)$ depending only on $n$ and $p$. We recall the following Proposition from [25].

**Proposition 3.8** Let $q \in L^\infty(\mathbb{S}^n)$, $r \in L^s(\mathbb{S}^n)$, for some $s > 1$, and let $1 < p \leq \frac{n+4}{n-4}$. Suppose $u \in H_2^2(\mathbb{S}^n)$ is a weak solution of the equation

(70) $$L_h^2 u = q |u|^{p-1} u + r \quad \text{on } \mathbb{S}^n.$$

Then for all $s > 1$, there exists a positive constant $\beta_{n,s}$ depending only on $n$ and $s$, such that if $\|q |u|^{2^\sharp - 2}\|_{L^{\frac{n}{4}}(\mathbb{S}^n)} \leq \beta_{n,s}$, then $u \in L^s(\mathbb{S}^n)$, and

$$\|u\|_s \leq C(n,s) \|r\|_s,$$

where $C(n,s)$ is a constant depending only on $n$ and $s$.

We are going to prove the following Proposition.

**Proposition 3.9** Let $k \in \mathbb{N}$, $a_1, \ldots, a_k \in \mathbb{S}^n$, $\alpha_1, \ldots, \alpha_k \in (0, +\infty)$, $A > 1, \gamma > 1$, and let $f \in C(\mathbb{S}^n)$ be a positive function. Suppose that

(71) $$A^{-1} \tau^\gamma \leq \lambda_i \leq A \tau^{\frac{1}{\gamma}} \quad i = 1, \ldots, k, \quad \tau \in (0, \tau_0),$$

(72) $$\alpha_i = \left(\frac{2}{n-4} d_n f(a^i) + o(1)\right)^{\frac{n-4}{2}} \quad \text{as } \tau \to 0,$$

and suppose $u$ is a solution of

$$P_h^n u = \frac{n-4}{2} f(x) |u|^{\frac{8}{n-4} - \tau} u \quad \text{on } \mathbb{S}^n, \quad \tau \in (0, \tau_0),$$

with

(73) $$u = \sum_{i=1}^k \alpha_i \delta_{a_i, \lambda_i} + v_\tau, \quad v_\tau \to 0 \text{ in } H_2^2(\mathbb{S}^n) \text{ as } \tau \to 0.$$

Then $u > 0$ for $\tau$ sufficiently small.

As an immediate consequence of Proposition 3.9 we have a complete proof of Proposition 3.6.

Now we come to the proof of Proposition 3.9. We are dealing with a solution $u$ of the equation

(74) $$P_h^n u = \frac{n-4}{2} f(x) |u|^{p-1} u \quad \text{on } \mathbb{S}^n,$$

where $p = \frac{n+4}{n-4} - \tau$. It is convenient to perform the conformal transformation $\varphi_{a, \lambda_1}$ on $\mathbb{S}^n$, which induces naturally the isometry $T_{\varphi_{a,\lambda_1}} : H_2^2(\mathbb{S}^n) \to H_2^2(\mathbb{S}^n)$ given by

$$T_{\varphi_{a,\lambda_1}} : u \quad \to \quad |\det \varphi_{a,\lambda_1}|^{\frac{n-4}{2n}} u \circ \varphi_{a,\lambda_1}.$$

Setting $\overline{u} = T_{\varphi_{a_1,\lambda_1}} u$, using (73) one can check that

$$\overline{u} = \alpha_1 + \sum_{i=2}^k \alpha_i \delta_{b_i, \zeta_i} + r_\tau,$$

where $\{b_i\}_i \subseteq \mathbb{S}^n$, $\zeta_i \to +\infty$ and $r_\tau \to 0$ in $H_2^2$ as $\tau \to 0$. Moreover, by the conformal invariance of $P_h^n$, $\overline{u}$ is a solution of

$$P_h^n \overline{u} = f(\varphi_{a_1,\lambda_1}(x)) \delta_{a_1,\lambda_1}^{-\tau} |\overline{u}|^{p-1} \overline{u} \quad \text{on } \mathbb{S}^n.$$

Now, writing $\overline{u} = \alpha_1 + w$, it is sufficient to prove that

(75) $$w \geq -\frac{1}{2}\alpha_1, \quad \text{for } \tau \text{ small enough.}$$

In fact this implies that $\overline{u} > 0$ and hence $u > 0$.
By a simple computation we obtain that $w$ satisfies

$$L_h^2 w = \phi_\tau(x, w(x)) \quad \text{in } \mathbb{S}^n,$$

where

$$\phi_\tau(x,t) = \frac{n-4}{2} f(\varphi_{a_1,\lambda_1}(x)) \, \delta_{a_1,\lambda_1}^{-\tau} \, |\alpha_1 + t|^{p-1}(\alpha_1 + t) + e_n\, t - d_n\, \alpha_1,$$

with $p = \frac{n+4}{n-4} - \tau$, and $e_n = \frac{c_n^2}{4} - d_n > 0$.
Setting $\tilde{f} = \phi_\tau(x, w(x))$, we denote by $w_1$ and $w_2$ the solutions of

$$L_h^2 w_1 = \tilde{f}^+, \quad L_h^2 w_2 = -\tilde{f}^-,$$

where $\tilde{f}^+ = \max\{\tilde{f}, 0\}$ and $\tilde{f}^- = -\min\{\tilde{f}, 0\}$. By the maximum principle, we have $w_1 \geq 0$, $w_2 \leq 0$; moreover, it is clear that $w_1 + w_2 = w$.
Inequality (75) is proved if we are able to show that

(76) $$w_2 \geq -\frac{1}{2}\alpha_1 \quad \text{for } \tau \text{ small.}$$

In order to do this, we set

$$\Xi = \{x \in \mathbb{S}^n : \tilde{f}(x) < 0\}.$$

We notice that $\phi_\tau(x, 0)$ is uniformly bounded on $\mathbb{S}^n$, and $\frac{\partial \phi}{\partial t}(x,t) \geq \gamma_1$ for a fixed $\gamma_1 > 0$, hence we have

(77) $$x \in \Xi \quad \Rightarrow \quad w(x) \leq C,$$

where $C$ is a fixed constant.

Fix a small $\varepsilon > 0$, and consider the sets

$$\Omega_\varepsilon = \Xi \cap \{x \in \mathbb{S}^n : -\varepsilon < w_2 \leq 0\}, \quad \Theta_\varepsilon = \Xi \cap \{x \in \mathbb{S}^n : w_2 \leq -\varepsilon\}.$$

**Lemma 3.10** *Let $\eta > 0$ be a small fixed positive number. The following properties hold true*

i) $\|w\|_{L^\infty(\Omega_\varepsilon)} \leq C$, *for some fixed constant $C$ ;*

ii) *in $\Theta_\varepsilon$, it is $\frac{|w|}{|w_2|} \leq C_\varepsilon$, for some constant $C_\varepsilon$ depending only on $\varepsilon$ ;*

iii) $\|w\|_{L^{\frac{2n}{n-4}}(\Xi)} \to 0$ *as $\tau \to 0$ ;*

iv) $|\tilde{f}| \leq C_1 + C_2 |w|^{\frac{n+4}{n-4}}$, *for some fixed positive constants $C_1$ and $C_2$ ;*

v) *for any $\eta > 0$, the function $\phi_\tau(x,t)$ satisfies the following properties*
   a) $\phi_\tau(x, 0) \to 0$ *uniformly on $\mathbb{S}^n \setminus B_\eta(-a_1)$ ,*
   b) $|\phi_\tau(x,t)| \leq |\phi_\tau(x,0)| + C\left(|t| + |t|^{\frac{n+4}{n-4}}\right)$, *for some fixed constant $C$ .*

PROOF  Property $i$) follows easily from (77) and from $w = w_1 + w_2 \geq w_2 \geq -\varepsilon > -1$ in $\Omega_\varepsilon$ (we can suppose $\varepsilon \in (0,1)$). Property $ii$) follows from the inequality $|w_2| \geq \varepsilon$ in $\Theta_\varepsilon$: in fact, in $\Xi$ we have $w_1 \leq C + w_2$, and hence we deduce immediately

$$\frac{|w|}{|w_2|} \leq \frac{|w_1| + |w_2|}{|w_2|} \leq \frac{C + 2|w_2|}{|w_2|} \leq 2 + \frac{C}{\varepsilon}.$$

Property $iii$) follows from (77) and $\zeta_i \to +\infty$. Properties $iv$) and $v$) are very easy to check, we just notice that for $v$)-$a$) we use $\delta^{-\tau}_{a_1, \lambda_1} \to 1$ uniformly on $\mathbb{S}^n$, because of (71) and (72). ∎

Before proving (76), we first show that $w_2$ tends to zero in $H^2_2(\mathbb{S}^n)$. From now on we write $B_\eta$ for $B_\eta(-a_1)$.

**Lemma 3.11** We have

$$\|w_2\|_{H^2_2} \to 0 \quad \text{as } \tau \to 0;$$

in particular, fixed $\varepsilon > 0$, $|\Theta_\varepsilon| \to 0$ as $\varepsilon \to 0$.

PROOF  We write

(78) $$L^2_h w_2 = \chi_\Xi \tilde{f} = \chi_\Xi \chi_{B_\eta} \tilde{f} + \chi_\Xi \chi_{\mathbb{S}^n \setminus B_\eta} \tilde{f}.$$

In $\Xi \cap B_\eta$ we can use property $iv$) above, so we deduce

$$\left\|\chi_\Xi \chi_{B_\eta} \tilde{f}\right\|_{L^{\frac{2n}{n+4}}} \leq \left(\int_{B_\eta \cap \Xi} \left(C_1 + C_2 |w|^{\frac{n+4}{n-4}}\right)^{\frac{2n}{n+4}}\right)^{\frac{n+4}{2n}} \leq C \left(|B_\eta| + \|w\|^{\frac{2n}{n-4}}_{L^{\frac{2n}{n-4}}(\Xi)}\right)^{\frac{n+4}{2n}}.$$

Since by $iii$) we have $\|w\|_{L^{\frac{2n}{n-4}}(\Xi)} \to 0$ as $\tau \to 0$, we obtain

(79) $$\left\|\chi_\Xi \chi_{B_\eta} \tilde{f}\right\|_{L^{\frac{2n}{n+4}}} = o_\eta(1) + o_\tau(1) \quad \text{for } \eta \text{ and } \tau \text{ small}.$$

We also have, by property $v$)-$b$)

$$\left|\chi_\Xi \chi_{\mathbb{S}^n \setminus B_\eta} \tilde{f}\right| \leq C\tau + C\left(|w| + |w|^{\frac{n+4}{n-4}}\right) \quad \text{in } \Xi \cap (\mathbb{S}^n \setminus B_\eta),$$

and hence

$$\left\|\chi_\Xi \chi_{\mathbb{S}^n \setminus B_\eta} \tilde{f}\right\|_{L^{\frac{2n}{n+4}}} \leq C \left(\int_{\mathbb{S}^n \setminus B_\eta} |\phi_\tau(x,0)|^{\frac{2n}{n+4}} + \int_\Xi \left(|w| + |w|^{\frac{2n}{n-4}}\right)\right)^{\frac{n+4}{2n}}.$$

Using $iii$) and $v$)-$a$) we deduce

(80) $$\left\|\chi_\Xi \chi_{\mathbb{S}^n \setminus B_\eta} \tilde{f}\right\|_{L^{\frac{2n}{n+4}}} \to 0 \quad \text{as } \tau \to 0.$$

From (79), (80) and the arbitrariness of $\eta$ it follows that $\lim_{\tau \to 0} \|\chi_\Xi \tilde{f}\|_{L^{\frac{2n}{n+4}}} = 0$. So the Lemma is a consequence of (78), (69) and the Sobolev embeddings. ∎

Now we come to the conclusion, namely we prove (76). We consider the function $\tilde{f}$ separately in the three sets $\Theta_\varepsilon, \Omega_\varepsilon \cap B_\eta$ and $\Omega_\varepsilon \cap (\mathbb{S}^n \setminus B_\eta)$.
In $\Theta_\varepsilon$ we have, using property $iv$)

$$|\tilde{f}(x)| \leq C_1 + C_2 |w|^{\frac{n+4}{n-4}} \leq C \left(\frac{1 + |w|^{\frac{n+4}{n-4}}}{|w_2|^{\frac{n+4}{n-4}}}\right) |w_2|^{\frac{n+4}{n-4}} \quad x \in \Theta_\varepsilon.$$

So, since $x \in \Theta_\varepsilon$, from $ii)$ it follows that

$$\tilde{f}(x) = g_{\varepsilon,\tau}(x) |w_2|^{\frac{8}{n-4}} w_2 \quad \text{in } \Theta_\varepsilon, \quad \text{and } |g_{\varepsilon,\tau}(x)| \leq C_\varepsilon, \tag{81}$$

where $C_\varepsilon$ is a positive constant depending only on $\varepsilon$.
By (77), $i)$ and $iv)$, we have

$$|\tilde{f}(x)| \leq C \quad x \in B_\eta \cap \Omega_\varepsilon, \tag{82}$$

for some fixed positive constant $C$.
Moreover for $x \in \Omega_\varepsilon$ it is $w(x) \geq \varepsilon$ hence, fixing $\eta > 0$, we have by $v) - a)$

$$-h^\eta_{\varepsilon,\tau} \leq \tilde{f}(x) < 0 \quad x \in (\mathbb{S}^n \setminus B_\eta) \cap \Omega_\varepsilon, \tag{83}$$

where $h^\eta_{\varepsilon,\tau}$ is a positive constant which tends to zero as $(\varepsilon, \tau) \to 0$.
Hence, taking into account (81), (82) and (83) we have

$$L^2_h w_2 = \chi_{\Theta_\varepsilon} g_{\varepsilon,\tau}(x) |w_2|^{\frac{8}{n-4}} w_2 + \hat{f}(x), \tag{84}$$

where

$$\|g_{\varepsilon,\tau}\|_{L^\infty} \leq C_\varepsilon, \quad |\hat{f}(x)| \leq h^\eta_{\varepsilon,\tau} + C \chi_{B_\eta}(x).$$

Now we fix $\alpha > \frac{n}{4}$ and $s > \frac{n+4}{n-4} \alpha$, so in particular it is $s > \alpha$. Since $w_2 \to 0$ in $H^2_2$ as $\tau \to 0$, by Lemma 3.11, it turns out that $\|\chi_{\Theta_\varepsilon} g_{\varepsilon,\tau}(x) |w_2|^{\frac{8}{n-4}}\|_{L^{\frac{n}{4}}} \to 0$ as $\tau \to 0$. We can apply Proposition 3.8 and we deduce that

$$\|w_2\|_s \leq C(n,s) \|\hat{f}\|_s. \tag{85}$$

The last $L^s$ norm can be estimated as

$$\|\hat{f}\|_s \leq \left( \int_{\mathbb{S}^n} \left( h^\eta_{\varepsilon,\tau} + C \chi_{B_\eta} \right)^s \right)^{\frac{1}{s}} \leq C \left( \|h^\eta_{\varepsilon,\tau}\|_\infty + \eta^{\frac{n}{s}} \right). \tag{86}$$

Applying (69) to equation (84) and using the Hölder inequality and the Sobolev embeddings we derive

$$\|w_2\|_\infty \leq C(n,\alpha) \left( \|\chi_{\Theta_\varepsilon} g_{\varepsilon,\tau} |w_2|^{\frac{n+4}{n-4}}\|_\alpha + C_{\alpha,s} \|\hat{f}\|_s \right). \tag{87}$$

Let $p = \frac{s}{\alpha} \frac{n-4}{n+4} > 1$, and let $p'$ be the conjugate exponent of $p$: using the Hölder inequality we have that

$$\|\chi_{\Theta_\varepsilon} g_{\varepsilon,\tau} |w_2|^{\frac{n+4}{n-4}}\|_\alpha \leq \left( \left( \int_{\mathbb{S}^n} |w_2|^s \right)^{\frac{1}{p}} \left( \int_\Xi |g_{\varepsilon,\tau}|^{\alpha p'} \right)^{\frac{1}{p'}} \right)^{\frac{1}{\alpha}} \leq C_\varepsilon^{\alpha p'} \|w_2\|_s^{\frac{s}{p\alpha}} |\Theta_\varepsilon|^{\frac{1}{\alpha p'}}. \tag{88}$$

From (87) and (88) it follows that

$$\|w_2\|_\infty \leq C(n,\alpha) C_\varepsilon^{\alpha p'} \|w_2\|_s^{\frac{s}{p\alpha}} |\Theta_\varepsilon|^{\frac{1}{\alpha p'}} + C(n,\alpha) C_{\alpha,s} \|\hat{f}\|_s.$$

Hence we have by (85) and (86)

$$\|w_2\|_\infty \leq C(n,s,\alpha) \left( C_\varepsilon^{\alpha p'} \left( \|h^\eta_{\varepsilon,\tau}\|_\infty + \eta^{\frac{n}{s}} \right)^{\frac{s}{p\alpha}} |\Theta_\varepsilon|^{\frac{1}{\alpha p'}} + \left( \|h^\eta_{\varepsilon,\tau}\|_\infty + \eta^{\frac{n}{s}} \right) \right).$$

Having fixed $\alpha$ and $s$, we can now choose first $\eta$ and then $\varepsilon$ such that $C(n,s,\alpha) \left( \|h^\eta_{\varepsilon,\tau}\|_\infty + \eta^{\frac{n}{s}} \right) \leq \frac{1}{4} \alpha_1$ for $\tau$ sufficiently small. Since $|\Theta_\varepsilon| \to 0$ as $\tau \to 0$ by Lemma 3.11, (76) follows from the last formula. This concludes the proof.

# 4  Appendix

## 4.1  A Pohozahev-type identity

**Proposition 4.1** *Let $n \geq 5$, let $B_r$ be the ball in $\mathbb{R}^n$ centered at $0$ and with radius $r$, and let $p \geq 1$ and let $u$ be a positive $C^4$ solution of*

$$\Delta^2 u = \frac{n-4}{2} f(x) u^q \quad x \in B_r. \tag{89}$$

*We have*

$$\frac{n-4}{2(q+1)} \sum_i \int_{B_r} x_i \frac{\partial f}{\partial x_i} u^{q+1}\, dx \;+\; \frac{n-4}{2}\left(\frac{n}{q+1} - \frac{n-4}{2}\right) \int_{B_r} f\, u^{q+1}\, dx$$

$$-\; r\frac{n-4}{2(q+1)} \int_{\partial B_r} f\, u^{q+1} d\sigma = \int_{\partial B_r} B(r, x, u, \nabla u, \nabla^2 u, \nabla^3 u)\, d\sigma,$$

*where*

$$B(r, x, u, \nabla u, \nabla^2 u, \nabla^3 u) \;=\; -\frac{n-2}{2}\Delta u \frac{\partial u}{\partial \nu} - \frac{r}{2}|\Delta u|^2 + \frac{n-4}{2} u \frac{\partial}{\partial \nu}(\Delta u)$$

$$+\; \langle x, \nabla u\rangle \frac{\partial}{\partial \nu}(\Delta u) - \Delta u \sum_i x_i \frac{\partial}{\partial \nu} u_i.$$

PROOF   Multiplying equation (89) by $u$ we have

$$\int_{B_r} u\, \Delta^2 u\, dx = \frac{n-4}{2}\int_{B_r} f\, u^{q+1}\, dx,$$

so integrating by parts we obtain

$$\int_{B_r} (\Delta u)^2\, dx = \frac{n-4}{2}\int_{B_r} f\, u^{q+1}\, dx - \int_{\partial B_r} \Delta u \frac{\partial u}{\partial \nu}\, d\sigma + \int_{\partial B_r} u \frac{\partial}{\partial \nu}(\Delta u)\, d\sigma.$$

Multiplying equation (89) by $\sum_{i=1}^n x_i u_i$ we obtain

$$\sum_{i=1}^n \int_{B_r} \Delta^2 u\, x_i\, u_i\, dx = \frac{n-4}{2} \sum_{i=1}^n \int_{B_r} x_i\, u_i\, f u^q\, dx. \tag{90}$$

Integrating by parts, we rewrite the right hand side of (90) as

$$\sum_{i=1}^n \int_{B_r} x_i u_i\, f u^q\, dx \;=\; -\frac{n}{q+1} \int_{B_r} f\, u^{q+1}\, dx + \frac{r}{q+1}\int_{\partial B_r} f\, u^{q+1}\, d\sigma$$

$$-\; \frac{1}{q+1} \sum_{i=1}^n \int_{B_r} x_i \frac{\partial f}{\partial x_i} u^{q+1}\, dx. \tag{91}$$

Next, one can transform the left hand side of (90) in the following way

$$\sum_{i=1}^n \int_{B_r} \Delta^2 u\, x_i\, u_i\, dx \;=\; \frac{4-n}{2}\int_{B_r}(\Delta u)^2\, dx + \frac{r}{2}\int_{\partial B_r}(\Delta u)^2\, d\sigma + \int_{\partial B_r} \Delta u \frac{\partial u}{\partial \nu}\, d\sigma \tag{92}$$

$$-\; \int_{\partial B_r} \langle x, \nabla u\rangle \frac{\partial}{\partial \nu}(\Delta u)\, dx + \sum_{i=1}^n \int_{\partial B_r} x_i \frac{\partial}{\partial \nu}(u_i)\, \Delta u\, dx.$$

So, putting together equations (90), (91) and (92) we obtain the result.  ∎

It is easy to check that the boundary term enjoys the following properties

**Corollary 4.2** $B(r, x, u, \nabla u, \nabla^2 u, \nabla^3 u)$ has the following properties

i) for $u(x) = |x|^{4-n}$, it is
$$B(r, x, u, \nabla u, \nabla^2 u, \nabla^3 u) = 0 \quad \text{for all } x \in \partial B_r,$$

ii) for $u(x) = |x|^{4-n} + A + \alpha(x)$, where $A > 0$ is some positive constant and $\alpha(x)$ is some $C^4$ function with $\alpha(0) = 0$, then there exists some $r^* > 0$ such that, for any $r$ with $0 < r < r^*$ we have
$$B(r, x, u, \nabla u, \nabla^2 u, \nabla^3 u) < 0 \text{ for all } x \in \partial B_r,$$

and
$$\lim_{r \to 0} \int_{\partial B_r} B(r, x, u, \nabla u, \nabla^2 u, \nabla^3 u) = -(n-4)^2 (n-2) \, \omega_{n-1} \, A.$$

## 4.2 A maximum principle for elliptic systems on domains

We recall the following result, see [36] page 193.

**Lemma 4.3** Let $D \in \mathbb{R}^n$ be a bounded smooth domain, and let $w(\cdot) > 0$ be a vector field on $\overline{D}$ such that
$$L_\mu[w_\mu] + \sum_{\nu=1}^{k} h_{\mu\nu} w_\nu \leq 0 \quad \text{in } D \quad \mu = 1, \ldots, k,$$

where
$$L_\mu = \sum_{i,j=1}^{n} a_{ij}^{(\mu)}(x) \frac{\partial^2}{\partial x_i \partial x_j} + \sum_{i=1}^{n} b_i^{(\mu)}(x) \frac{\partial}{\partial x_i},$$

where $a_{ij}^{(\mu)}(x)$, $b_i^{(\mu)}(x)$ are uniformly bounded, $a_{ij}^{(\mu)}(x)$ are uniformly elliptic, and $h_{\mu\nu} \geq 0$ for $\mu \neq \nu$. Suppose that the vector field $z(x)$ satisfies the system of inequalities
$$L_\mu[z_\mu] + \sum_{\nu=1}^{k} h_{\mu\nu} z_\nu \geq 0 \quad \text{in } D \quad \mu = 1, \ldots, k,$$

and that there exists some constant $M > 0$ such that $z \leq M w$ on $\partial D$. Then $z \leq M w$ in $\overline{D}$.

## 4.3 Some properties of biharmonic functions

We recall the following well known Lemma, see for example [7].

**Lemma 4.4** (Böcher) Suppose $n \geq 3$, $a_0 \in \mathbb{R}$, and $v \in C^2(B_2 \setminus \{0\})$ satisfies the conditions
$$\begin{cases} \Delta v = 0 & \text{in } B_2 \setminus \{0\}, \\ v \geq a_0 |x|^{2-n} & \text{in } B_2 \setminus \{0\}. \end{cases}$$

Then there exist $a_1 \geq a_0$ and an harmonic function $d : B_1 \to \mathbb{R}$ such that
$$v(x) = a_1 |x|^{2-n} + d(x) \quad x \in B_1 \setminus \{0\}.$$

Taking into account Lemma 4.4, we can prove the following analogous result regarding the biharmonic operator.

**Lemma 4.5** *Suppose $n \geq 5$, and suppose $v \in C^4(B_1 \setminus \{0\})$ satisfies the conditions*

$$\begin{cases} \Delta^2 v = 0 & \text{in } B_2 \setminus \{0\}, \\ v \geq 0 & \text{in } B_2 \setminus \{0\}, \\ \Delta v \geq 0 & \text{in } B_2 \setminus \{0\}. \end{cases}$$

*Then there exist $a_1, a_2 \geq 0$ a function $b \in C^\infty(B_1)$ with $\Delta^2 b = 0$ such that*

$$v(x) = a_1 |x|^{4-n} + a_2 |x|^{2-n} + b(x) \quad x \in B_{\frac{1}{2}} \setminus \{0\}.$$

PROOF  Set $w = \Delta v$. Then one can easily check that $w$ satisfies the assumptions of Lemma 4.4 with $a_0 = 0$, so there exist $a_3 \geq 0$ and a function $d \in C^\infty(B_1)$ with $\Delta d = 0$ such that

$$w(x) = a_3 |x|^{2-n} + d(x) \quad x \in B_1 \setminus \{0\}.$$

Define $\tilde{v} : B_1 \setminus \{0\}$ to be

$$\tilde{v}(x) = \frac{a_3}{2(n-4)} |x|^{4-n} + \Delta^{-1} d(x) \quad x \in B_1 \setminus \{0\}.$$

where $\Delta^{-1} d$ denotes a classical solution of $\Delta u = d$ in $B_1$. It is easy to see that, setting $\overline{v} = v - \tilde{v}$, there exists $\overline{C} > 0$ such that $\overline{v}$ satisfies

$$\begin{cases} \Delta \overline{v} = 0 & \text{in } B_1 \setminus \{0\}, \\ \overline{v} \geq -\overline{C} |x|^{4-n} & \text{in } B_1 \setminus \{0\}. \end{cases}$$

Hence, applying Lemma 4.4 to $\overline{v}$, there exist $a_2 \geq 0$ and $e \in C^\infty(B_1)$ with $\Delta e = 0$ such that

$$\overline{v}(x) = a_2 |x|^{2-n} + e(x) \quad \text{in } B_1 \setminus \{0\}.$$

Hence the Lemma follows setting

$$a_1 = \frac{a_3}{2(n-4)}, \quad b = e + \Delta^{-1} d.$$

The proof is concluded. ■

# References

[1] Andersson S.I., Doebner H.D. (Eds), Nonlinear partial differential operators and quantization procedures. Proceedings of a workshop held at the Technische Universitt Clausthal, Clausthal, July 1981. Edited by S. I. Andersson and H.-D. Doebner.

[2] Ambrosetti A., Garcia Azorero J., Peral A., Perturbation of $-\Delta u + u^{\frac{(N+2)}{(N-2)}} = 0$, the Scalar Curvature Problem in $\mathbb{R}^N$ and related topics, Journal of Functional Analysis, 165 (1999), 117-149.

[3] Arapostathis A., Ghosh M.F., Marcus S., Harnack's inequality for cooperative weakly coupled systems, Communications in Partial Differential Equations, 24, 9-10 (1999), 1555-1571.

[4] Aubin T., Meilleures constantes dans le théorème d' inclusion de Sobolev et un théorème de Fredholm non linéaire pour la transformation conforme de la courbure scalaire, Journal of Functional Analysis, 32, 1979, 148-174.

[5] Aubin T., Some Nonlinear Problems in Differential Geometry, Springer-Verlag, 1998.

[6] Aubin T. and Bahri A. Méthodes de topologie algébrique pour le probléme de la courbure scalaire prescrite, Journal des Mathématiques Pures et Appliquées, 76 (1997), 525-549.

[7] Axler S., Bourdon P., Ramey W., *Harmonic Function Theory*, Springer-Verlag, GTM 137, 1992.

[8] Bahri A., Critical points at infinity in some variational problems, Research Notes in Mathematics, 182, Longman-Pitman, London, 1989.

[9] Bahri A. An invariant for Yamabe type flows with applications to scalar curvature problems in higher dimensions, Duke Mathematical Journal, 81 (1996), 323-466.

[10] Bahri A., Coron J.M., The Scalar-Curvature problem on the standard three-dimensional sphere, Journal of Functional Analysis, 95 (1991), 106-172.

[11] Ben Ayed M., Chen Y., Chtioui H., Hammami M., On the prescribed scalar curvature problem on 4-manifolds, Duke Mathematical Journal, 84 (1996), 633-677.

[12] Branson T.P., Group representations arising from Lorentz conformal geometry, Journal of Functional Analysis, 74 (1987), 199-291.

[13] Branson T.P., Differential operators canonically associated to a conformal structure, Mathematica Scandinavica, 57-2 (1985), 293-345.

[14] Branson T.P., Chang S.Y.A., Yang P.C., Estimates and extremal problems for the log-determinant on 4-manifolds, Communications in Mathematical Physics, 149 (1992), 241-262.

[15] Chang S.Y.A., On a fourth order PDE in conformal geometry, Preprint, 1997.

[16] Chang S.Y.A., Gursky M.J., Yang P.C., The scalar curvature equation on 2- and 3- spheres, Calculus of Variations and Partial Differential Equations, 1 (1993), 205-229.

[17] Chang S.Y.A., Gursky M.J., Yang P.C., Regularity of a fourth order non-linear PDE with critical exponent, American Journal of Mathematics, 121-2 (1999), 215-257.

[18] Chang S.Y.A., Qing J., Yang P.C., On the Chern-Gauss-Bonnet integral for conformal metrics on $\mathbb{R}^4$, Duke Mathematical Journal, to appear.

[19] Chang S.Y.A., Qing J., Yang P.C., Compactification for a class of conformally flat 4-manifolds, Inventiones Mathematicae, to appear.

[20] Chang S-Y.A., Yang P., A perturbation result in prescribing scalar curvature on $\mathbb{S}^n$, Duke Mathematical Journal, 64 (1991), 27-69.

[21] Chang S.Y.A., Yang P.C., Extremal metrics of zeta functional determinants on 4-manifolds, Annals of Mathematics, 142, 1995, 171-212.

[22] Chang S.Y.A., Yang P.C., On a fourth order curvature invariant, Comtemporary Mathematics, 237, Spectral Problems in Geometry and Arithmetic, Ed. T. Branson, AMS, 1999, 9-28.

[23] Connes A., Noncommutative geometry. Academic Press, Inc., San Diego, CA, 1994.

[24] Djadli Z., Hebey E., Ledoux M., Paneitz-type operators and applications, Duke Mathematical Journal, 104-1 (2000) 129-169.

[25] Djadli Z., Malchiodi A., Ould Ahmedou M., *Prescribed fourth order conformal invariant on the standard sphere, Part I*, preprint.

[26] Doebner H.D., Hennig J.D. (Eds), Differential geometric methods in mathematical physics. Proceedings of the twelfth international conference held at the Technical University of Clausthal, Clausthal, August 30–September 2, 1983. Edited by H.-D. Doebner and J. D. Hennig.

[27] Escobar J., Schoen,R., Conformal metrics with prescribed scalar curvature, Inventiones Mathematicae, 86 (1986), 243-254.

[28] Gilbarg D., Trudinger N., Elliptic Partial Differential Equations of Second Order, 2nd edition, Springer-Verlag, 1984.

[29] Gursky M.J., The Weyl functional, de Rham cohomology, and Kahler-Einstein metrics, Annals of Mathematics, 148, 1998, 315-337.

[30] Li Y.Y., Prescribing scalar curvature on $\mathbb{S}^n$ and related topics, Part I, Journal of Differential Equations, 120 (1995), 319-410; Part II, Existence and compactness, Communications in Pure and Applied Mathematics, 49 (1996), 437-477.

[31] Lin C.S., A classification of solutions of conformally invariant fourth order equation in $\mathbb{R}^n$, Commentari Mathematici Helveticii, 73 (1998), 206-231.

[32] Hebey E., Changements de métriques conformes sur la sphère - Le problème de Nirenberg, Bull. Sci. Math. 114 (1990), 215-242.

[33] Moser J., On a nonlinear problem in differential geometry, Dynamical Systems (M. Peixoto ed.), Academic Press, New York, 1973, 273-280.

[34] Paneitz S., A quartic conformally covariant differential operator for arbitrary pseudo-Riemannian manifolds, Preprint, 1983.

[35] Paneitz S., Essential unitarization of symplectics and applications to field quantization, Journal of Functional Analysis, 48-3, 1982, 310–359.

[36] Protter M.H., Weinberger H.F., Maximum Principles in Differential Equations, Springer-Verlag, 2nd edition, 1984.

[37] Schoen R., On the number of constant scalar curvature metrics in a conformal class, in "Differential Geometry: A Symposium in Honor of Manfredo Do Carmo", (H.B.Lawson and K.Tenenblat Editors), (1991), 331-320, Wiley, New York.

[38] Schoen R., Zhang D., Prescribed scalar curvature on the $n$-sphere, Calculus of Variations and Partial Differential Equations, 4 (1996), 1-25.

[39] Serrin J., Zou H., Non-existence of positive solutions of Lane-Emden systems, Differential and Integral Equations, 9-4 (1996), 635-653.

[40] Wei J., Xu X.: On conformal deformations of metrics on $\mathbb{S}^n$, Journal of Functional Analysis, 157-1 (1998), 292-325.

**Z.D.** : Université de Cergy-Pontoise - Département de Mathématique - Site de Saint-Martin - 2 Avenue Adolphe Chauvin - F 95302 Cergy-Pontoise Cedex, France
**A.M.** : Rutgers University - Department of Mathematics - Hill Center - Busch campus - 110 Frelinghuysen Rd - Piscataway, NJ08854-8019 , USA
**M.O.A.** : Scuola Internazionale Superiore di Studi avanzati (SISSA) - Via Beirut, 2-4 - 34014 Trieste - Italy

e-mail : Zindine.Djadli@math.u-cergy.fr   ,   malchiod@math.rutgers.edu   ,   ahmedou@sissa.it